\documentclass[11pt,A4]{article}
\usepackage{subfigure}
\usepackage{graphicx,psfrag,epsfig,epsf}
\usepackage{amsmath}
\usepackage{amsfonts}
\usepackage{amssymb}
\usepackage{colordvi}
\usepackage{cite}
\newtheorem{lemma}{Lemma}
\newtheorem{definition}{Definition}
\newtheorem{theorem}{Theorem}
\newtheorem{proposition}{Proposition}

\newenvironment{proof}{\paragraph{Proof:}}{\medskip\endpf}
\def\endpf{\proofbox\medskip\\}
\newcommand{\proofbox}{\hfill\mbox{\rule{2mm}{2mm}}}
\newcommand{\reals}{I\!\!R}
\newcommand{\natnums}{I\!\!N}
\newcommand{\bb}{(\cdot)}
\newcommand{\eqbyd}{\triangleq}
\def\veps{\varepsilon}

\def\bmbf#1{\bar{\mathbf{#1}}}
\def\mbb#1{\mathbb{#1}}
\def\mc{\mathcal}
\def\mbf#1{\mathbf{#1}}
\def\tmbf#1{\tilde{\mathbf{#1}}}

\def\tends{\rightarrow}
\def\interior#1{\mathrm{interior}\{{#1}\}}

\def\co{\mathrm{co}}

\def\int#1{\overset{\circ}{#1}}

\def\sstar#1{\overset{\ast\ }{#1}}
\def\proj{\mathrm{Proj}}
\numberwithin{equation}{section}
\begin{document}
\title{Characterization of the solution to a constrained $H_\infty$
  optimal control problem\thanks{Research supported by Engineering
    and Physical Sciences Research Council, UK and the Royal Academy
  of Engineering.}\\[1ex] \small{
Technical Report CUED/F-INFENG/TR.518,\\ Department of
Engineeering, University of Cambridge}}

\author{D. Q.  Mayne\thanks{Department of Electrical and Electronic
    Engineering, Imperial College London.}, S. V. Rakovi\'{c}$^\dag$, R.
    B. Vinter$^\dag$ and E. C. Kerrigan\thanks{Royal Academy of
    Engineering Post-doctoral Research Fellow, Department of
    Engineering, University of Cambridge.} }

\date{30 March 2004\\Revised 14 February  2005}

\maketitle

\begin{abstract}
  This paper characterizes the solution to a finite horizon min-max
  optimal control problem where the system is linear and discrete-time
  with control and state constraints, and the cost quadratic; the
  disturbance is negatively costed, as in the standard $H_\infty$
  problem, and is constrained. The cost is minimized over
  control policies and maximized over disturbance sequences so that
  the solution yields a feedback control. It is shown that the value
  function is piecewise quadratic and the optimal control policy
  piecewise affine, being quadratic and affine, respectively, in
  polytopes that partition the domain of the value function.
\end{abstract}
\normalsize

\section{Introduction}
Characterizations of solutions to constrained optimal control problems
appeared in the papers
\cite{dona:goodwin:1999,seron:goodwin:dona:2000,bemporad:morari:dua:pistikopoulos:2002,mayne:rakovic:2003a}
that deal with the constrained linear-quadratic problem, in the papers
\cite{bemporad:borrelli:morari:2000a,bemporad:borrelli:morari:2000c,mayne:rakovic:2003a,mayne:rakovic:2002b}
and thesis \cite{borrelli:2002} that deal with hybrid or piecewise
affine systems, and in papers that deal with min-max optimal control
problems
\cite{ramirez:camacho:2001,kerrigan:mayne:2002a,bemporad:borrelli:morari:2003,kerrigan:maciejowski:2004,diehl:bjornberg:2004}. In these
papers it is shown that the value function is piecewise affine or
piecewise quadratic (depending on the nature of the cost function in
the optimal control problem) and the control law is piecewise affine,
being quadratic or affine in polytopes that constitute a polytopic
partition of the domain of the value function. When disturbances are
present, it is necessary to compute the solution sequentially using
dynamic programming as in \cite{kerrigan:mayne:2002a}. In this paper,
which is motivated by recent research on $H_\infty$ model predictive
control
\cite{chen:scherer:allgower:1997,magni:nicolao:scattolini:1998,magni:nijmeijer:vanderschaft:2001,lee:kouvaritakis:2000a,magni:nicolao:scattolini:allgower:2003,grimm:teel:zacharian:2003,lee:2003,kim:2004},
we obtain an explicit characterization of the solution to a
constrained, min-max optimal control problem and consider here the
choice of terminal cost and constraint set to ensure stability of the
closed loop system with receding horizon control.  The term $H_\infty$
is used somewhat loosely since we consider the min-max problem with
fixed $\gamma$. We consider, therefore, the problem of controlling a
linear, discrete-time system described by
\begin{equation}
  \label{eq:1.1}
  x^+=Ax+Bu+Gw, \qquad y=Cx+Du
\end{equation}
where $x \in \reals^n$ is the state, $u \in \reals^m$ the control and
$w \in \reals^p$ an additive disturbance (the `adversary'); $x^+$ is
the successor state and $y \in \reals^r$ is the costed output. We
frequently write the system dynamics in \eqref{eq:1.1} in the form
\begin{equation*}
  x^+=f(x,u,w)
\end{equation*}
where $f(x,u,w)\eqbyd Ax+Bu+Gw$.  The system is subject to hard
control and state constraints
\begin{equation}
  \label{eq:1.2}
  u \in U, \qquad x \in X
\end{equation}
where $U\subseteq \reals^m$ is a (compact) polytope and $X\subseteq
\reals^n$ a polytope; each set contains the origin in its interior
(the assumption that $X$ is a polytope rather than a
polyhedron\footnote{A polyhedron is a set described by a finite set of
  inequalities; a polytope is a bounded polyhedron.} is made for
simplicity).  The disturbance $w$ is constrained to lie in the
polytope $W\subseteq \reals^p$; $W$ contains the origin in its
interior.

Let $\pi \eqbyd \{\mu_0\bb,\mu_1\bb,\ldots,\mu_{N-1}\bb\}$ denote a
control policy (sequence of control \emph{laws}) over horizon $N$ and
let $\mbf{w} \eqbyd \{w_0,w_1,\ldots,w_{N-1}\}$ denote a sequence of
disturbances. Also, let $\phi(i;x,\pi,\mbf{w})$ denote the solution of
\eqref{eq:1.1} when the initial state is $x$ at time $0$, the control
policy is $\pi$ and the disturbance sequence is $\mbf{w}$, so that
$\phi(i;x,\pi,\mbf{w})$ is the solution, at time $i$ of
\begin{eqnarray}
  \label{eq:1.3}
  x_{i+1} &=& Ax_i + B \mu_i(x_i)+Gw_i\\
  \label{eq:1.4}
  x_0 &=& x
\end{eqnarray}
The cost $V_N(x,\pi,\mbf{w})$, if the initial state is $x$, the control
policy $\pi$ and the disturbance sequence $\mbf{w}$, is
\begin{equation}
  \label{eq:1.5}
  V_N(x,\pi,\mbf{w}) \eqbyd \sum_{i=0}^{N-1} \ell(x_i,u_i,w_i) + V_f(x_N)
\end{equation}
where, for all $i$, $x_i\eqbyd \phi(i;x,\pi,\mbf{w})$ and $u_i\eqbyd
\mu_i(x_i)$; $V_f\bb$ is a terminal cost that may be chosen, together
with a terminal constraint set $X_f$ defined below, to ensure
stability of the resultant receding horizon controller (see \S6).  The
stage cost $\ell\bb$ is a quadratic function, positive definite in $x$
and $u$, and negative definite in $w$:
\begin{equation}
  \label{eq:1.6}
\ell(x,u,w) \eqbyd (1/2)|x|_Q^2+(1/2)|u|_R^2 -(\gamma^2/2) |w|^2
\end{equation}
where $\gamma >0$, $|z|_Z^2 \eqbyd z'Zz$,  and $Q$ and $R$ are
positive definite. The stage cost may be expressed as
\begin{equation}
  \label{eq:1.7}
\ell(x,u,w) \eqbyd (1/2)|y|^2 -(\gamma^2/2) |w|^2,\ y \eqbyd Hz
\end{equation}
where $z\eqbyd(x,u)$ and $H$ is a suitably chosen matrix ($(x,u)$
should be interpreted as a column vector $(x',u')'$ in matrix
expressions).  The terminal cost $V_f\bb$ is a quadratic function
\begin{equation}
  \label{eq:1.8}
V_f(x) \eqbyd(1/2) |x|_{P_f}^2
\end{equation}
in which $P_f$ is positive definite.  The optimal control problem
$\mathbb{P}_N(x)$ that we consider is
\begin{equation}
 \label{eq:1.9}
 \mathbb{P}_N(x):\qquad V_N^0(x) = \inf_{\pi \in \Pi_N(x)} \max_{\mbf{w} \in \mc{W}}
 V_N(x,\pi,\mbf{w})
\end{equation}
where $\mc{W} \eqbyd W^N$, is the set of admissible disturbance
sequences, and $\Pi_N(x)$ is the set of admissible policies, i.e.
those policies that satisfy, for all $\mbf{w} \in \mc{W}\eqbyd W^N$,
the state and control constraints \eqref{eq:1.2}, and the terminal
constraint
\begin{equation}
  \label{eq:1.10}
   x_N \in X_f.
\end{equation}
Inclusion of the hard disturbance constraint $\mbf{w} \in \mc{W}$
is necessary when state constraints are present since,
otherwise, for any policy $\pi$ chosen by the controller, we can
expect that there exists a disturbance sequence $\mbf{w}$ that
transgresses the state constraint.  The terminal constraint set is a
polytope, containing the origin in its interior, that satisfies $X_f
\subseteq X$, ensuring satisfaction of the state constraint at time
$N$. Hence the set of admissible policies is
\begin{multline}
  \label{eq:1.11}
  \Pi_N(x) \eqbyd \{\pi \mid \phi(i;x,\pi,\mbf{w}) \in X,\
  \mu_i(\phi(i;x,\pi,\mbf{w})) \in U, i=0,1,\ldots,N-1,\\
  \phi(N;x,\pi,\mbf{w}) \in X_f,\ \forall \mbf{w} \in \mc{W}\}
\end{multline}
Let $X_N$ denote the set of initial states for which a solution to
$\mathbb{P}_N(x)$ exists (the domain of $V_N^0\bb$, the
controllability set), i.e.
\begin{equation}
  \label{eq:1.12}
  X_N \eqbyd \{x \mid \Pi_N(x) \not = \emptyset\}.
\end{equation}
In addition to characterizing the solution to a min-max optimal
control problem that has not previously been characterized, this paper
provides an improvement of the transformation procedure used in
\cite{mayne:rakovic:2003a,mayne:rakovic:2002b,borrelli:2002} to obtain
a parametric solution to the optimal control problem; the improvement
simplifies the determination of the polytopes in which the control law
and value function are affine and quadratic respectively and avoids
unnecessary sub-partitioning of overlapping polytopes required in
\cite{mayne:rakovic:2003a,borrelli:2002}.

\section{Dynamic Programming for Constrained Problems}
The solution to $\mathbb{P}_N(x)$ may be obtained as follows. For all
$j \in \natnums_+ \eqbyd\{1,2,\ldots\}$, let problem $\mathbb{P}_j$,
the partial return function $V_j^0\bb$, and the controllability set
$X_j$ be defined as in \eqref{eq:1.5}--\eqref{eq:1.9} with $j$
replacing $N$; $j$ denotes ``time-to-go''. Then the sequences
$\{V_j^0\bb, \kappa_j\bb, X_j\}$, where $\kappa_j\bb$ denotes
the optimal control law $\mu_{N-j}^0\bb$ at time
$i=N-j$, may be calculated recursively as follows
\cite{mayne:2001a,kerrigan:mayne:2002a}:
\begin{alignat}{1}
  \label{eq:2.1}
  V_j^0(x) &= \min_{u \in U} \max_{w \in W}\{\ell(x,u,w)+
           V_{j-1}^0(f(x,u,w)) \mid f(x,u,W) \subseteq X_{j-1}\}\\
  \kappa_j(x) &= \arg\min_{u \in U}  \max_{w \in
  W}\{\ell(x,u,w)+V_{j-1}^0(f(x,u,w)) \mid \nonumber\\
  \label{eq:2.2}
  & \qquad \qquad \qquad  \qquad \qquad \qquad \qquad \qquad  \qquad
   f(x,u,W) \subseteq X_{j-1}\}\\
  \label{eq:2.3}
   X_j &=X \cap\{x\mid \exists u \in U\ \mathrm{such\ that}\
  f(x,u,W) \subseteq X_{j-1}\}
\end{alignat}
with boundary conditions
\begin{equation}
  \label{eq:2.4}
  V_0^0(x) = V_f(x), \qquad X_0 = X_f.
\end{equation}
The condition $f(x,u,W) \subseteq X_{j-1}$ in \eqref{eq:2.2} and
\eqref{eq:2.3} may be expressed as
\begin{equation}
  \label{eq:2.5}
   Ax+Bu \in X_{j-1} \ominus GW
\end{equation}
where $\ominus$ denotes Pontryagin set difference defined by
$\mc{A}\ominus \mc{B} \eqbyd \{x \mid \{x\} \oplus \mc{B} \subseteq
\mc{A}\}$ ($\oplus$ denotes set addition). For each integer $j$ let
$Z_j \subseteq \reals^n\times\reals^m$ be defined by
\begin{equation}
  \label{eq:2.6}
 Z_j \eqbyd \{(x,u) \in X \times U \mid f(x,u,W) \subseteq X_{j}\}
\end{equation}
so that, from \eqref{eq:2.3},
\[
X_j = \proj _X Z_{j-1}.
\]
Here, and in the sequel, if a set $Z$, say, lies in a product space
$\reals^n \times \reals^m$, $\proj_X:\reals^n \times \reals^m \tends
\reals^n$ denotes the projection operator defined by $\proj_X Z =
\{x \mid \exists u \in \reals^m \text{\ such\ that\ } (x,u) \in Z\}$
($\reals^n$ is regarded as $x$-space). Similarly, if $\Phi$ is a set
in the product space $\reals^n \times \reals^m \times \reals ^p$,
$\proj_Z \Phi$ denotes the set $\{{z} \mid \exists w \in \reals^p
\text{\ such\ that\ } (z,w) \in \Phi\}$.  We can now establish some
preliminary properties of the solution to $\mathbb{P}_N$. To analyze
$\mathbb{P}_N(x)$ it is convenient to introduce the functions
$J_j^0\bb$, $j=1,2,\ldots$, defined by
\begin{equation}
  \label{eq:2.7}
  J_j^0(x,u) \eqbyd \max_{w \in W} \{\ell(x,u,w)+V_{j}^0(f(x,u,w))\}
\end{equation}
The recursive equations \eqref{eq:2.1}-\eqref{eq:2.3} may therefore be
rewritten as
\begin{alignat}{1}
  \label{eq:2.8}
  V_j^0(x) &= \min_{u \in U} \{J_{j-1}^0(x,u) \mid f(x,u,W) \subseteq
  X_{j-1}\}\\
  \label{eq:2.9}
  J_{j-1}^0(x,u) &\eqbyd  \max_{w \in W}
  \{\ell(x,u,w)+V_{j-1}^0(f(x,u,w))\}\\
  \label{eq:2.10}
  \kappa_j(x) &= \arg\min_{u \in U} \{J_{j-1}^0(x,u) \mid f(x,u,W)
  \subseteq X_{j-1}\}\\
  \label{eq:2.11}
  X_j &=X \cap\{x\mid \exists u \in U\ \mathrm{st}\ f(x,u,W) \subseteq
  X_{j-1}\}
\end{alignat}
for $j=1,\ldots,N$ with endpoint conditions $V_0\bb=V_f\bb$,
$X_0=X_f$.  Under our assumptions the sets $X_j$ and $Z_j$ are
compact. If $X_0=X_f$ is robust control invariant, the sets $X_j$ are
nested ($X_j \supseteq X_{j-1}$ for all $j \ge 1$). For each $j$, the
domain of $V_j^0\bb$ includes $X_j$ but we are only interested in its
values on $X_j$; similarly, the domain of $J_j^0\bb$ includes $Z_j$
but we are only interested in its values on $Z_j$.

\section{Parametric Optimization}
We seek a parametric solution to problem $\mathbb{P}_N(x)$, i.e. a
solution for all values of the parameter which, in this case, is the
state $x$. More precisely, since we employ constrained dynamic
programming, we seek a parametric solution to problems
$\mathbb{P}_j(x)$ for all $j \in \{1,\ldots,N\}$.  First, we introduce
a few useful definitions.
\begin{definition}
  \label{d:1}
For any positive integer $J$, $\mc{I}_J \eqbyd \{1,2,\ldots,J\}$; for
any set $\mc{X}$, $\mc{J}^\mc{X}$ denotes an index set associated with
a partition of $\mc{X}$.
\end{definition}
\begin{definition}
  \label{d:2}
  A set $\mc{P}=\{P_i \mid i \in \mc{J}\}$, for some index set
  $\mc{J}$, is called a polyhedral (polytopic) partition of a closed
  (compact) set $\mc{X}$ if $\mc{X}= \cup_{i \in \mc{J}} P_i$, and the
  sets $P_i,\ i \in \mc{J}$ are polyhedrons (polytopes) with non-empty
  interiors which are non-intersecting ($\mathrm{interior}({P}_i) \cap
  \mathrm{interior}({P}_j)=\emptyset$ for all $i,j \in \mc{J}, i\neq
  j$).
\end{definition}
\begin{definition}
  \label{d:3}
  A function $V:\mc{X} \tends \reals$ is said to be continuous
  piecewise quadratic on a polyhedral (polytopic) partition
  $\mc{P}=\{P_i \mid i \in \mc{J}\}$ of $\mc{X}$ if it is
  \emph{continuous} and satisfies
\[
V(x) =(1/2)|x|_{Q_i}^2+q_i'x+r_i, \qquad \forall x \in P_i,\  i
\in \mc{J}
\]
for some $Q_i, q_i, r_i$,  $i \in \mc{J}$.
Similarly, a function $\kappa:\mc{X} \tends \mc{U}$ is said to be
piecewise affine on a polyhedral partition $\mc{P}=\{P_i \mid i \in
\mc{J}\}$  of $\mc{X}$ if it is continuous and satisfies
\[
\kappa(x) = K_ix+k_i, \qquad \forall x \in P_i,\  i
\in \mc{J},
\]
for some $K_i, k_i$,  $i \in \mc{J}$, where $\mc{P}$ has the
properties specified above.
\end{definition}

The dynamic programming recursion \eqref{eq:2.8}-\eqref{eq:2.11}
requires the repeated solution of two prototype problems
$\mathbb{P}_\mathrm{min}$ and $\mathbb{P}_\mathrm{max}$ defined next:
\begin{alignat}{1}
  \label{eq:3.1}
  \mathbb{P}_\mathrm{min}(x): \qquad V^0(x) &= \min_{u} \{J(x,u)
  \mid (x,u) \in \mc{Z}\}\\
  \label{eq:3.2}
  \mathbb{P}_\mathrm{max}(z): \qquad J^0(z) &= \max_{w} \{V(z,w)  \mid
  w \in W\}
\end{alignat}
The minimizer in $\mathbb{P}_\mathrm{min}(x)$ and the maximizer in
$\mathbb{P}_\mathrm{max}(z)$ are defined, respectively, by
\begin{alignat}{1}
 \label{eq:3.3}
 \kappa(x) & \eqbyd \arg \min_{u} \{J(x,u) \mid (x,u) \in \mc{Z}\}\\
 \label{eq:3.4}
 \nu(z) & \eqbyd \arg \max_{w} \{V(z,w)  \mid  w \in W\}.
\end{alignat}
Problem $\mathbb{P}_\mathrm{min}(x)$ is the prototype for Problem
\eqref{eq:2.8} with $V^0(x)$ replacing $V_j^0(x)$, $J(x,u)$ replacing
$J_{j-1}^0(x,u)$, and $(x,u) \in \mc{Z}$ replacing the constraints
$f(x,u,W) \subseteq X_{j-1}$ ($Ax+Bu \in X_{j-1}\ominus GW$)
\emph{and} $u \in U$.  Similarly Problem $\mathbb{P}_\mathrm{max}(z)$
is the prototype for Problem \eqref{eq:2.9} with $J^0(z)$ replacing
$J_{j-1}^0(z)$, $V(z,w)$ replacing $\ell(x,u,w)+V_{j-1}^0(f(x,u,w))$,
and $z$ replacing $(x,u)$.  We first obtain the parametric solution of
$\mathbb{P}_\mathrm{min}$.

\subsection{The minimization problem $\mathbb{P}_\mathrm{min}$}

\begin{figure}[htbp]
   \psfrag{Z}{$\mc{Z}$}
   \psfrag{G1}{$G_1$}
   \psfrag{G2}{$G_2$}
   \psfrag{x}[t][]{$x$}
   \psfrag{u}[r][]{$u$}
   \psfrag{U}[l][]{$\mc{U}(x_1)$}
   \psfrag{z1}[][]{$\ z_1$}
   \psfrag{z2}[][]{$z_2$}
   \psfrag{x1}[][]{$B(x_1,\veps_1') \cap\mc{X}$}
   \psfrag{bx}[][r]{$B(x_2,\veps_2')\cap\mc{X}\ \ $}
   \begin{center}
   \includegraphics[width=8cm]{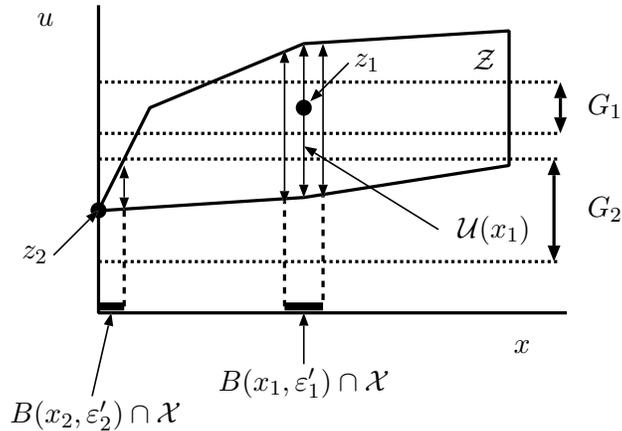}
   \caption{Inner  semi-continuity of $\mc{U}(x)$}
   \label{f:1}
   \end{center}
 \end{figure}

The solution to $\mathbb{P}_\mathrm{min}(x)$ has properties given in
Proposition \ref{p:1} that has a simpler hypothesis than previous
versions of this result. For completeness, continuity of the control
law is also proven.
\begin{proposition}
  \label{p:1}
  Suppose $J: \mc{Z} \tends \reals$ is a strictly convex, continuous
  function and that $\mc{Z}$ is a polytope.  Then, for all $x \in
  \mc{X}=\mathrm{Proj}_X \mc{Z}$, the solution $\kappa(x)$ to
  $\mathbb{P}_\mathrm{min}(x)$ exists and is unique.  The value
  function $V^0\bb$ is strictly convex and continuous with domain
  $\mc{X}$, and the control law $\kappa\bb$ is continuous on $\mc{X}$.
\end{proposition}
\begin{proof}
  For all $x \in \mc{X}$, $\mc{U}(x) \eqbyd \{u \mid (x,u)
  \in \mc{Z}\}$ is convex and compact.  Let $\Lambda :=
  \{(\lambda_1,\lambda_2) \mid \lambda_1 \ge 0, \lambda_2 \ge 0,
  \lambda_1+\lambda_2=1\}$. For all $x_1,x_2$ in $\mc{X}$, all
  $\lambda=(\lambda_1,\lambda_2) \in \Lambda$:
\begin{eqnarray*}
 V^0(\lambda_1 x_1+\lambda_2 x_2) &=&\min_u \{J(\lambda_1 x_1+\lambda_2
x_2,u) \mid (\lambda_1 x_1+\lambda_2 x_2,u) \in \mc{Z}\} \\
&\le& J(\lambda_1 x_1+\lambda_2 x_2,\lambda_1 u_1+\lambda_2 u_2\},\
 u_i \eqbyd  \kappa(x_i),\ i=1,2\\
&=& J(\lambda_1(x_1,u_1) +\lambda_2(x_2,u_2))
\end{eqnarray*}
But $\lambda_1(x_1,u_1) +\lambda_2(x_2,u_2) \in \mc{Z}$ since $\mc{Z}$ is
convex and $(x_i,u_i)\in \mc{Z}$, $i=1,2$.  Since $J\bb$ is strictly
convex
\begin{eqnarray*}
V^0(\lambda_1 x_1+\lambda_2 x_2) &\le&
\lambda_1 J(x_1,u_1) + \lambda_2 J(x_2,u_2)\\
&=& \lambda_1 V^0(x_1) + \lambda_2 V^0(x_2)\ \forall \lambda_1,\ \lambda_2
\in \Lambda
\end{eqnarray*}
where the last inequality is strict if $\lambda_1 \not \in
\{0,1\}$ so that $V^0\bb$ is strictly convex. Since $J\bb$ is
strictly convex, $\kappa(x)$ is unique at each $x \in \proj_X
\mc{Z}$.

The constraint $(x,u) \in \mc{Z}$ imposes an implicit state-dependent
constraint $u \in \mc{U}(x)$ on $u$ where the set-valued function
$\mc{U}\bb$ is defined by
\begin{equation*}
  \mc{U}(x) \eqbyd \{u \mid (x,u) \in \mc{Z}\}.
\end{equation*}
We claim that $\mc{U}\bb$ is continuous (both outer and inner
semi-continuous on $\mc{X}=\proj_X\mc{Z}$, the domain of $\mc{U}\bb$.
By definition \cite{polak:1997}, the set-valued map $\mc{U}\bb$ is
outer semi-continuous at $x \in \mc{X}$ if $\mc{U}(x)$ is closed and
if, for any compact set $G$ such that $\mc{U}(x) \cap G = \emptyset$
there exists an $\varepsilon >0$ such that $\mc{U}(x) \cap G =
\emptyset$ for all $x' \in B(x,\veps)\cap \mc{X}$. The set-valued map
$\mc{U}\bb$ is inner semi-continuous at $x \in \mc{X}$ if, for any
open set $G\subseteq \reals^m$ such that $G\cap \mc{U}(x)\neq
\emptyset$, there exists an $\varepsilon >0$ such that $G\cap
\mc{U}(x') \neq \emptyset$ for all $x' \in B(x,\veps)\cap \mc{X}$.
Here $B^j(x,\veps) \eqbyd \{x'\in \reals^j\mid |x'-x| \le \veps\}$.
The set-valued map $\mc{U}\bb$ is outer semi-continuous because its
graph, $\mc{Z}$, is closed so that, given any sequence $\{(x_i,u_i)\}$
in $\mc{Z}$ ($u_i \in \mc{U}(x_i)$ for all $i$) such that $(x_i,u_i)
\tends (\bar{x},\bar{u})$, we have $(\bar{x},\bar{u}) \in \mc{Z}$ so
that $\bar{u} \in \mc{U}(\bar{x})$. Hence $\mc{U}\bb$ is outer
semi-continuous \cite{polak:1997}.  We can establish inner
semi-continuity using the following result \cite{clarke:2005} whose
proof is given in the appendix.

\begin{lemma}\emph{(Clarke).}
   \label{l:1}
   Suppose $\mc{Z}$ is a polytope in $\reals^n \times \reals^m$ and
   let $\mc{X}$ denote its projection on $\reals^n$ ($\mc{X}= \{x \mid
   \exists u \in \reals^m \textrm{\ such\ that\ } (x,u) \in
   \mc{Z}\}$).  Let $\mc{U}(x) \eqbyd \{u \mid (x,u) \in \mc{Z}\}$.
   Then there exists a $K > 0$ such that, for all $x, x' \in \mc{X}$,
   for all $u \in \mc{U}(x)$, there exists a $u' \in \mc{U}(x')$ such
   that $|u'-u| \le K|x'-x|$.
\end{lemma}
Let $x, x'$ be arbitrary points in $\mc{X}$ and $\mc{U}(x)$ and
$\mc{U}(x')$ the associated sets (Figure \ref{f:1} illustrates the
proof for two cases: $x=x_1$ and $x=x_2$). Let $G$ be an open set such
that $\mc{U}(x) \cap G \neq \emptyset$ and let $u$ be an arbitrary
point in $\mc{U}(x) \cap G$. Because $G$ is open, there exist an
$\veps>0$ such that $B(u,\veps) \eqbyd \{v \mid |v-u| \leq \veps\}
\subset G$. Let $\veps' \eqbyd \veps/K$. From Lemma \ref{p:1}, there
exists a $u' \in \mc{U}(x') \cap G$ for all $x' \in
B(x,\veps')\cap\mc{X}$. This implies $\mc{U}(x') \cap G \neq
\emptyset$ for all $x' \in B(x,\veps')\cap\mc{X}$, so that $\mc{U}\bb$
is inner semi-continuous.
\end{proof}

To solve the parametric problem $\mathbb{P}_\mathrm{min}$, we
develop further the reverse transformation procedures proposed in
\cite{mayne:2001c} and utilized in
\cite{mayne:rakovic:2002b,borrelli:2002,mayne:rakovic:2003a}. We
assume that $J\bb$ is strictly convex and continuous piecewise
quadratic on a polytopic partition $\mc{P}^\mc{Z} =\{P_i^\mc{Z}
\mid i \in \mc{J}^\mc{Z}\}$ (for some index set $\mc{J}^\mc{Z}$)
of $\mc{Z}$. For each $z=(x,u) \in \mc{Z}$, let $S^\mc{Z}(z)$, the
index set of active polytopes at $z$, be defined by
\begin{equation}
  \label{eq:3.5}
  S^\mc{Z}(z) \eqbyd \{i \in \mc{J}^\mc{Z} \mid z \in P_i^\mc{Z}\}
\end{equation}
so that $S^\mc{Z}(z)$ is the set of indices of active polytopes at
$z$.  Similarly, for each $x \in \mc{X }\eqbyd \proj_X(\mc{Z})$, let
$S_{\mc{Z}}^0(x)$ be defined by
\begin{equation}
  \label{eq:3.6}
  S_{\mc{Z}}^0(x) \eqbyd S^\mc{Z}(x,\kappa(x))
\end{equation}
where $\kappa(x)$, the solution of $\mathbb{P}_{\min}(x)$, is defined by
\begin{equation}
  \label{eq:3.7}
  \kappa(x) \eqbyd \arg \min_{u} \{J(x,u) \mid (x,u) \in \mc{Z}\}
\end{equation}
so that $S_\mc{Z}^0(x)$ is the index set of polytopes active at the
solution $\kappa(x)$ of $\mathbb{P}_\text{min}(x)$. For each $i \in
\mc{J}^\mc{Z}$ we consider the simpler problem $\mathbb{P}_i(x)$
defined by
\begin{alignat}{1}
  \label{eq:3.8}
  \mathbb{P}_i(x): \qquad V_i^0(x) &= \min_u \{J_i(x,u) \mid
  (x,u) \in  P_i^\mc{Z}\}\\
  \label{eq:3.9}
  \kappa_i(x) &= \arg \min_u \{J_i(x,u) \mid (x,u) \in  P_i^\mc{Z}\}
\end{alignat}
where $J(z)=J_i(z)$ on $P_i^\mc{Z}$ and
\begin{equation}
  \label{eq:3.10}
 J_i(z) =  (1/2)z'Q_iz+q_i'z+s_i
\end{equation}
for some $Q_i, q_i, s_i$, all $z=(x,u) \in P_i^\mc{Z}$.  For each $i$,
problem $\mathbb{P}_i(x)$ is a quadratic program since $J\bb$
is quadratic on the polytope $P_i^\mc{Z}$.
\begin{proposition}
  \label{p:2}
  Suppose $J: \mc{Z} \tends \reals$ is strictly convex and continuous
  and that $\mc{Z}$ is a polytope with a polytopic partition
  $\{P_i^\mc{Z} \mid i \in \mc{J}^\mc{Z}\}$.  Then, for all $x \in
  P_i^\mc{X}\eqbyd\mathrm{Proj}_X {P}_i^\mc{Z}$, all $i \in
  \mc{J}^\mc{Z}$, the solution $\kappa_i(x)$ to $\mathbb{P}_i(x)$
  exists and is unique. Moreover, the value function $V_i^0\bb$ is
  strictly convex and continuous, and $\kappa_i\bb$ is continuous, in
  $P_i^\mc{X}$.
\end{proposition}
The proof of Proposition \ref{p:2} is almost identical to the proof of
Proposition \ref{p:1} noting that, for each $i \in \mc{J}^\mc{Z}$,
$P_i^\mc{Z}$ is a polytope.  Problem $\mathbb{P}_i(x)$, although
simpler, is an artificial problem since the constraint $(x,u) \in
P_i^\mc{Z}$ does not appear in the original problem
$\mathbb{P}_\text{min}(x)$, so that it is not obvious how best to
relate the solutions to problems $\mathbb{P}_i(x)$, $i \in
\mc{J}^\mc{Z}$, to the solution of $\mathbb{P}_\text{min}(x)$.  This
difficulty was not totally satisfactorily dealt with in the literature
quoted above but is resolved in the following result.
\begin{proposition}
  \label{p:3}
  Suppose $J: \mc{Z} \tends \reals$ is continuous piecewise quadratic
  and strictly convex on a polytopic partition $\mc{P}^\mc{Z}$ of
  $\mc{Z}$.  Then $u$ is optimal for the minimization problem
  $\mathbb{P}_\mathrm{min}(x)$ if and only if $u$ is optimal for the
  problems $\mathbb{P}_i(x)$ (i.e. if and only if
  $\kappa(x)=\kappa_i(x)$) for all $i \in {S}_{\mc{Z}}^0(x)$.
\end{proposition}
\begin{proof}
  Suppose $u=\kappa(x)$ is optimal for $\mathbb{P}_\text{min}(x)$ but
  that, contrary to what is to be proven, there exists an $i \in
  \mc{S}_{\mc{Z}}^0(x)$ such that $u$ is not optimal for
  $\mathbb{P}_i(x)$.  Let $u_i=\kappa_i(x)$ denote the solution of
  $\mathbb{P}_i(x)$. By definition, $(x,u_i) \in P_i^\mc{Z}$ and
  $(x,u) \in P_i^\mc{Z}$ (since $(x,u) \in P_j^\mc{Z}$ for all $j \in
  \mc{S}_{\mc{Z}}^0(x)$).  Hence $u_i =\kappa_i(x)$ satisfies
  $V_i^0(x) =J_i(x,u_i)=J(x,u_i) < J_i(x,u) = J(x,u)=V^0(x)$ where we
  have made use of the fact that $J(x,v) = J_i(x,v)$ if $(x,v) \in
  P_i^\mc{Z}$. Hence $J(x,u_i) < V^0(x)$ which contradicts the
  optimality of $u$ for $\mathbb{P}_\text{min}(x)$.  Suppose, next,
  that $u$ is optimal for $\mathbb{P}_i(x)$ for all $i \in
  \mc{S}_{\mc{Z}}^0(x)$ (so that $\kappa_i(x)=u$ for all $i \in
  \mc{S}_{\mc{Z}}^0(x)$) but that, contrary to what is to be proved,
  $u$ is not optimal for $\mathbb{P}_\text{min}(x)$ so that there
  there exists a $u^\ast$ satisfying $(x,u^\ast) \in \mc{Z}$ and
  $J(x,u^\ast) < J(x,u)$.  Because $u \in P_i^\mc{Z}$ for all $i\in
  \mc{S}_{\mc{Z}}^0(x)$ and $d(u,P_j^\mc{Z}) >0$ for all $j\in
  \mc{J}^\mc{Z} \setminus \mc{S}_{\mc{Z}}^0(x)$, there exists a
  $\lambda \in (0,1]$ and an $i \in \mc{S}_{\mc{Z}}^0(x)$ such that
  $u_\lambda \eqbyd u +\lambda(u^\ast-u)$ satisfies $(x,u_\lambda) \in
  P_i^\mc{Z}$.  Since $u \mapsto J(x,u)$ is convex and $J(x,u^\ast) <
  J(x,u)$ it follows that $J(x,u_\lambda) < J(x,u)$.  But
  $J(x,u_\lambda) = J_i(x,u_\lambda)$ (since $(x,u_\lambda) \in
  P_i^\mc{Z}$) and $J(x,u)= J_i(x,u)$ (since $(x,u) \in P_i^\mc{Z}$)
  so that $J_i(x,u_\lambda) <J_i(x,u)$, a contradiction of the
  optimality of $u$ for $\mathbb{P}_i(x)$ for all $i \in
  \mc{S}_{\mc{Z}}^0(x)$.
\end{proof}
Proposition \ref{p:3} shows that the solution to
$\mathbb{P}_\text{min}(x)$ is also the solution to a set of quadratic
programs, namely $\mathbb{P}_i(x)$ for $i \in \mc{S}_{\mc{Z}}^0(x)$.
We now analyse problem $\mathbb{P}_i(x)$ in more detail. Suppose that,
for each $i \in \mc{J}^\mc{Z}$, polytope $P_i^\mc{Z}$ is defined by
\begin{equation}
  \label{eq:3.11}
  P_i^\mc{Z} \eqbyd \{z=(x,u) \mid M_iu \le N_ix+p_i\}
\end{equation}
where $M_i,N_i,p_i$ each have $r_i$ rows, so that
\begin{alignat}{1}
 \nonumber
 \mathbb{P}_i(x): \qquad &V_i^0(x) = \min_u \{J_i(x,u) \mid M_iu
 \le N_ix+p_i\}\\
  \label{eq:3.12}
  & \kappa_i(x) = \arg\min_u \{J_i(x,u) \mid M_iu \le N_ix+p_i\}.
\end{alignat}
The $j^\mathrm{th}$ constraint  $M_i^ju \le N_i^jx+p_i^j$ is said
to be active at $(x,u)$ if $ M_i^ju = N_i^jx+p_i^j$. The set of
active constraints for  $\mathbb{P}_i(x)$ is $I_i^0(x)$, the set of
constraints active at $(x,\kappa_i(x))$, so that
\begin{equation}
  \label{eq:3.18}
  I_i^0(x) \eqbyd \{ j \in \mc{I}_{r_i} \mid M_i^j\kappa_i(x) =
  N_i^jx+p_i^j\}.
\end{equation}
where the superscript $j$ on a matrix (or vector) denotes the
$j^\text{th}$ row of the matrix (or vector). It follows from the
definition of $\kappa_{i,I}(x) = K_{i,I}x+k_{i,I}$ that $I_i^0(x) = I$
for all $x \in \interior{X}_{i,I}$ and that $I_i^0(x) \subseteq I$ on
the boundary of $X_{i,I}$.
The solution to $\mathbb{P}_i(x)$ is simple if the set of active
constraints $I_i^0(x)$ for the problem is known in advance
\cite{seron:dona:goodwin:2000,bemporad:morari:dua:pistikopoulos:2002}.
Suppose therefore the set of active constraints for $\mathbb{P}_i(x)$
at $(x,\kappa_i(x))$ ($\kappa_i(x)$ is the solution of
$\mathbb{P}_i(x)$) is known, apriori, to be $I$, i.e. $I^0(x)=I$.  Then
$ \mathbb{P}_i(x)$ is replaced by the simpler, equality constrained,
problem
\begin{equation}
  \label{eq:3.13}
  \mathbb{P}_{i,I}(x): \qquad V_{i,I}^0(x) = \min_u \{J_i(x,u) \mid
  M^j_iu = N_i^jx+p_i^j,\ j \in I \}
\end{equation}
This is a quadratic optimization problem with affine equality
constraints; the solution to this problem has, as is well known, the
form
\begin{alignat}{1}
  \label{eq:3.14}
  V_{i,I}^0(x) &=(1/2)x'Q_{i,I}x+q_{i,I}'x+s_{i,I}\\
  \label{eq:3.15}
  \kappa_{i,I}(x) &= K_{i,I}x+k_{i,I}
\end{alignat}
for some $Q_{i,I}$, $s_{i,I}$, $K_{i,I}$ and $k_{i,I}$.  Let $M_i^I$
denote the matrix with rows $M_i^j$, $j \in I$.  Let $PC_{i,I} \eqbyd
\{(M_i^I)'\lambda \mid \lambda \ge 0\} \subseteq \reals^m$ denote the
polar cone at $0$ to the cone $\mc{F}_{i,I} = \{h \mid M^j_i h \le 0,\
j \in I\}$ of feasible directions $h$ for problem $\mc{P}_i(x)$ at
$u=\kappa_{i,I}(x)$. The polar cone depends \emph{solely} on $I$, the
set of active constraints; it does not depend on the parameter $x$;
also \cite{mayne:rakovic:2003a} $-\nabla_u J_i(x,\kappa_{i,I}(x)) \in
PC_{i,I}$ if and only if $\langle\nabla_u J_i(x,\kappa_{i,I}(x)),
h\rangle \ge 0$ for all feasible directions $h$ ($h \in
\mc{F}_{i,I}$), i.e.  if and only if $ \kappa_{i,I}(x)$ is optimal for
the problem $\min_u \{J_i(x,u) \mid M^j_iu = N_i^jx+p_i^j,\ j \in I,
M_i^ju \le N_i^jx+p_i^j,\ j \in \mc{I}_{r_i} \setminus I\}$.  For each
$i \in \mc{J}^\mc{Z}$, let $P_i^\mc{X}$ denote the polytope defined by
\begin{equation}
  \label{eq:3.16}
  P_i^\mc{X} \eqbyd \{x \mid \exists u \text{ s.t. } (x,u) \in P_i^\mc{Z}\}
  = \proj_X(P_i^\mc{Z})
\end{equation}
The polytope $P_i^\mc{X}$ is the domain of $V_i^0\bb$.
The following result holds \cite{mayne:rakovic:2003a}:
\begin{proposition}
  \label{p:4}
  The affine control law $\kappa_{i,I}\bb$ is optimal for problem
  $\mathbb{P}_i(x)$, at all $x$ in the polytope $X_{i,I}\subseteq
  P_i^\mc{X}$  defined by
\begin{equation}
  \label{eq:3.17}
X_{i,I} \eqbyd \left\{ x \in X \ \pmb{\Bigg\vert} \begin{array}{rcl}
  M_i^j(K_{i,I}x+k_{i,I}) &\le& N_i^jx+p_i^j,\ j \in
  \mc{I}_{r_i} \setminus I\\
   -\nabla_u J_i(x,\kappa_{i,I}(x))&\in&PC_{i,I} \end{array} \right\}
\end{equation}
\end{proposition}

The restriction $x \in X$ is included in the definition of $X_{i,I}$
since it is not included as a constraint in $\mathbb{P}_i(x)$). Since
the affine control law $\kappa_{i,I}\bb$ is such that the equality
constraint $M_i \kappa_{i,I}(x) = N_ix+p_i$ is satisfied for all $x$,
and since the last inequality in \eqref{eq:3.17} ensures that
$\kappa_{i,I}(x)$ is optimal for $\min_u \{J_i(x,u) \mid M^j_iu =
N_i^jx+p_i^j,\ j \in I, M_i^ju \le N_i^jx+p_i^j,\ j \in \mc{I}_{r_i}
\setminus I\}$, it follows that $\kappa_{i,I}(x)$ is optimal for
$\mathbb{P}_i(x)$ in the polytope $X_{i,I}$.  Thus
\cite{seron:dona:goodwin:2000,bemporad:morari:dua:pistikopoulos:2002,mayne:rakovic:2003a}
the solution to $\mathbb{P}_i(x)$ is affine, and the value function
quadratic, in each polytope $X_{i,I}$; the set of all such non-empty
polytopes (as $I$ ranges over the subsets of $\{1,2,\ldots,r_i\}$
constitute a polytopic partition of $\proj_X(P_i^\mc{Z})$ so that the
solution $\kappa_i\bb$ to $\mathbb{P}_i(x)$ is piecewise affine, and
the value function $V_i^0\bb$ is piecewise quadratic, on this
polytopic partition.

However, in our case, since we have to `marry' a set of polytopes
$X_{i,I_i}$ for all $i$ such that polytope $P_i^\mc{Z}$ is active
(i.e. $i \in \mc{S}_\mc{Z}^0(x)$), is active, it is preferable to
parameterize the polytopes in which the solution to $\mathbb{P}_i(x)$
is affine by the state $\bar{x}$, say, at which $I_i$ is active rather
than by the set $I_i$ of active constraints.  Also, for each $\bar{x}
\in \mc{X}=\proj_X(\mc{Z})$, let the polytope $X_i(\bar{x}) \subseteq
P_i^\mc{X}$ be defined by
\begin{equation}
  \label{eq:3.19}
  X_i(\bar{x})\eqbyd   X_{i,I_i^0(\bar{x})}
\end{equation}
where, for each set $I\subseteq \mc{J}^\mc{Z}$, $X_{i,I}$ is defined
by \eqref{eq:3.17}.  It follows from \eqref{eq:3.17}, with $I$ replaced
by $I_i^0(\bar{x})$, that $\bar{x} \in  X_i(\bar{x})$. It was shown
above that $I_i^0(x)=I_i^0(\bar{x})$ for all $x$ in the interior of
$X_i(\bar{x})$; it follows from \eqref{eq:3.19} that
$X_i(x)=X_i(\bar{x})$ for all $x$ in the interior of $X_i(\bar{x})$.

\begin{figure}[htbp]
  \psfrag{P1}{$P_1^\mc{Z}$}
  \psfrag{P2}{$P_2^\mc{Z}$}
  \psfrag{z}[t][]{$x$}
  \psfrag{w}[r][]{$u$}
  \psfrag{m}[r][]{$\kappa(x)$}
  \psfrag{m1}[r][]{$\kappa_1(x)$}
  \psfrag{m2}[r][]{$\kappa_2(x)$}
  \psfrag{W}{$U$}
  \psfrag{z11}[b][l]{$X_{11}$}
  \psfrag{z12}[b][l]{$X_{12}$}
  \psfrag{z13}[b][]{$X_{13}$}
  \psfrag{z21}[t][]{$X_{21}$}
  \psfrag{z22}[t][]{$X_{22}$}
  \psfrag{z23}[t][]{$X_{23}$}
  \subfigure[Solutions of $\mathbb{P}_i(x)$, $i=1,2$]
  {\includegraphics[width=5cm]{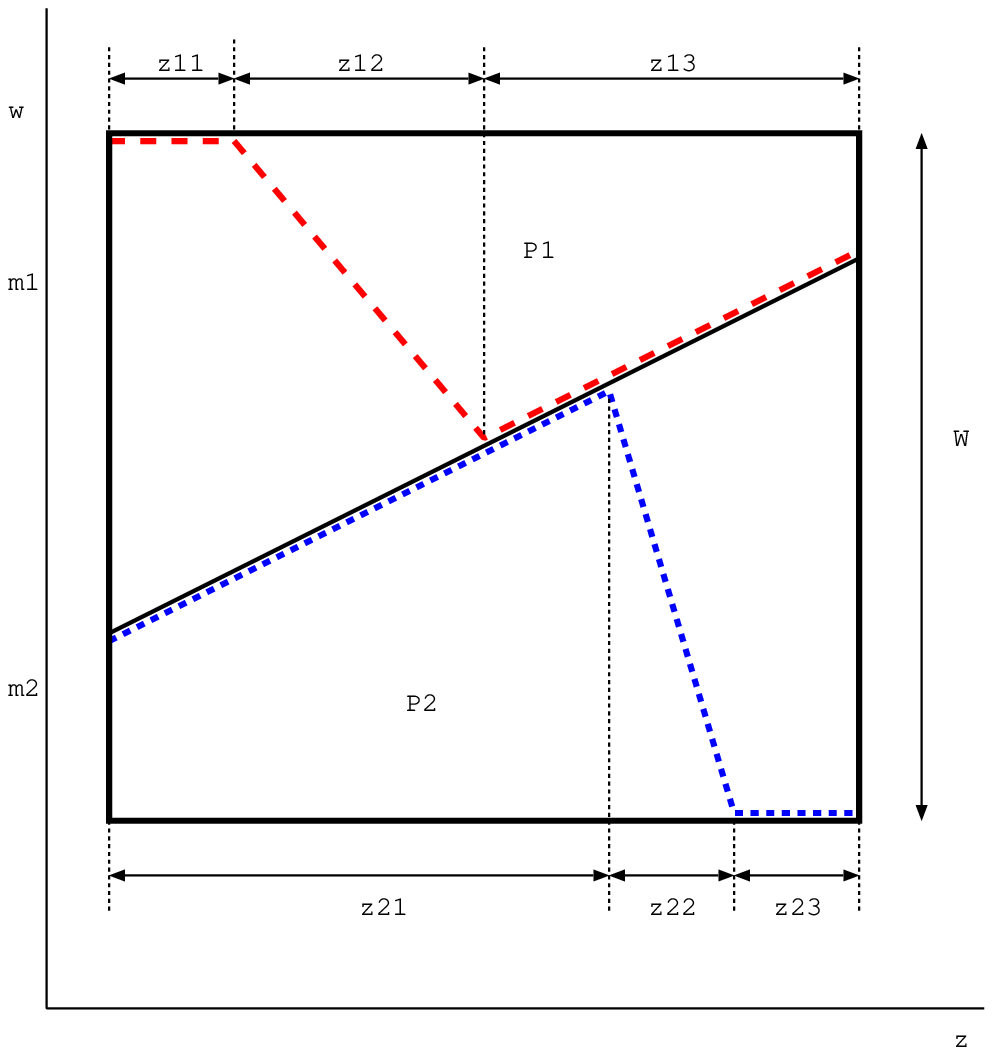}}
  \hspace{1in}
  \subfigure[Solution of $\mathbb{P}_\text{min}(x)$]
  {\includegraphics[width=5cm]{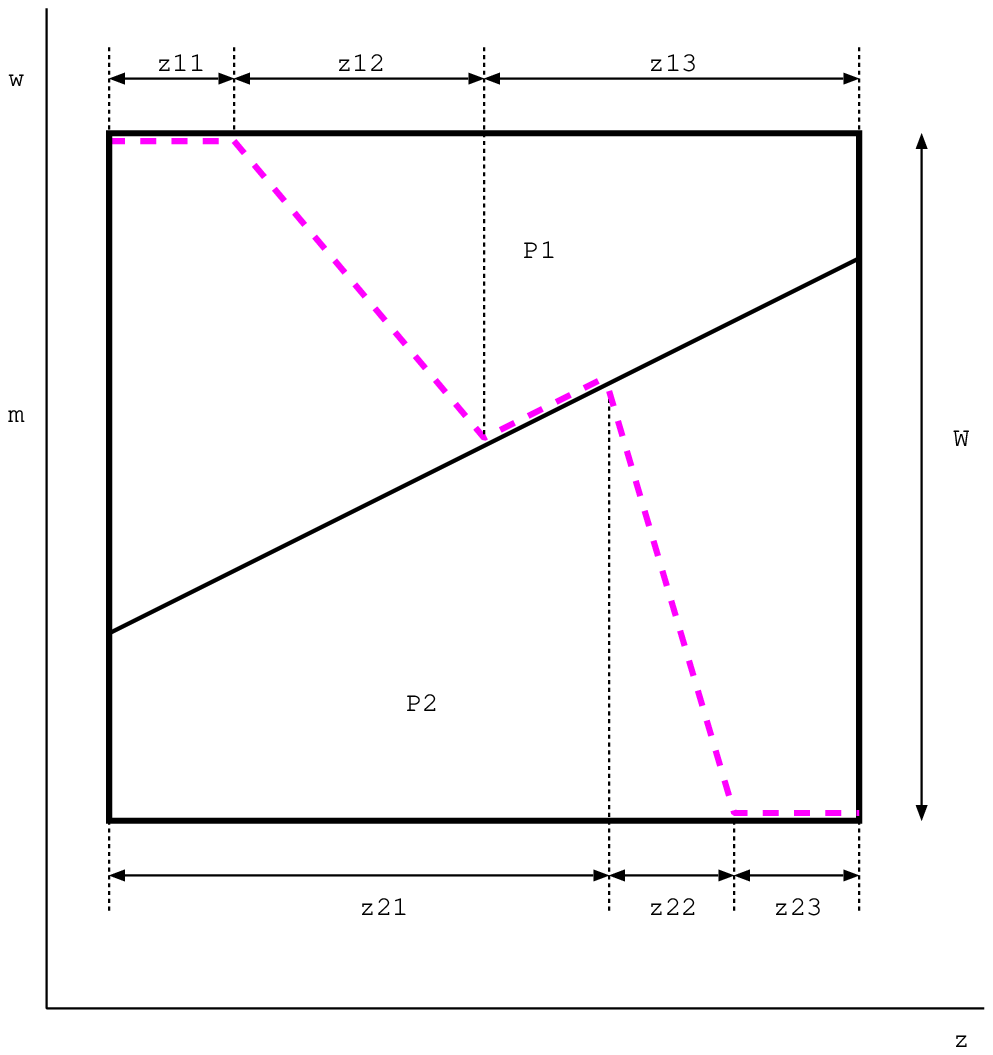}}
  \caption{Problem $\mathbb{P}(x)$}
  \label{f:2}
\end{figure}

The polytope $X(\bar{x})$ that figures in the
parametric solution of $\mathbb{P}_\text{min}(x)$ is defined, for each
$\bar{x}\in \mc{X}$, by
\begin{equation}
  \label{eq:3.20}
  X(\bar{x}) \eqbyd \bigcap \{X_i(\bar{x}) \mid i \in
  S_{\mc{Z}}^0(\bar{x})\}
\end{equation}
Clearly $\bar{x} \in X(\bar{x})$ and $X(x')=X(\bar{x})$ for all $x'
\in X(\bar{x})$. For each $\bar{x} \in \mc{X}$, let the functions
$V_{\bar{x}}^0(\cdot)$ and $\kappa_{\bar{x}}(\cdot)$ be
defined on $X(\bar{x})$ by
\begin{alignat}{1}
  \label{eq:3.21}
  V_{\bar{x}}^0(x) &\eqbyd V_i^0(x),\ \forall i \in
  S_{\mc{Z}}^0(\bar{x})\\
  \label{eq:3.22}  \kappa_{\bar{x}}(x) &\eqbyd \kappa_i(x),\ \forall i \in
  S_{\mc{Z}}^0(\bar{x})
\end{alignat}
The domain of each function is $X(\bar{x})$; that the functions are
well defined follows from Proposition \ref{p:3} and equations
\eqref{eq:3.19} and \eqref{eq:3.20} which show that
$\kappa_i(x)=\kappa_j(x)$ for all $i,j \in S_{\mc{Z}}^0(\bar{x})$, all
$x \in X(\bar{x})$.  Summarizing, we have:
\begin{theorem}
  \label{th:1}
  Suppose $J:{\mc{Z}} \tends \reals$ is continuous piecewise quadratic and
  strictly convex on a polytopic partition $\mc{P}^\mc{Z}=\{P_i^\mc{Z}
  \mid i \in \mc{J}^\mc{Z}\}$ of ${\mc{Z}}$ and that, for each $i \in
  \mc{J}^\mc{Z}$, $P_i^\mc{Z}$ has a non-empty interior. Then the
  value function $V^0\bb$ is continuous piecewise quadratic and strictly convex
  on a polytopic partition $\mc{P}^\mc{X}=\{P_i^\mc{X} \mid i \in
  \mc{J}^\mc{X}\}$ of $\mc{X}$$=\proj_X(\mc{Z})$. The minimizer
  $\kappa\bb$ is piecewise affine on $\mc{P}^\mc{X}$. The polytopes
  $P_i^\mc{X}$ are each of the form $X(\bar{x})$ for some $\bar{x} \in
  \mc{X}$; the value function and optimal control law satisfy $V^0(x)
  = V_{\bar{x}}^0(x)$ and $\kappa(x)=\kappa_{\bar{x}}(x)$ for all $x
  \in X(\bar{x})$ and some $\bar{x} \in \mc{X}$.
\end{theorem}
The proof of this result follows from Propositions \ref{p:3} and
\ref{p:4} and the discussion above. The result is illustrated in
Figure \ref{f:2} for the simple case when $\mc{Z}=
\mc{P}_1^\mc{Z}\cup\mc{P}_2^\mc{Z}$ has two partitions
$\mc{P}_1^\mc{Z}$ and $\mc{P}_2^\mc{Z}$ in each of which $J\bb$ is
quadratic. Problem $\mc{P}_1^\mc{Z}(x)$ is, therefore, a parametric
quadratic program; its solution $\kappa_1(x)$ is known to be piecewise
affine on a polytopic partition of $\mc{P}_1^\mc{X}$; in Figure
\ref{f:2}, $\mc{P}_1^\mc{X}=\mc{X}$ and the polytopic partition is
$\{X_{11}, X_{12},X_{13}\}$. The solution to $\mc{P}_1^\mc{Z}(x)$ is
$\kappa_1(x)$ which is affine in each of the polytopes $X_{11},
X_{12}$ and $X_{13}$. Similarly the solution to the quadratic program
$\mc{P}_1^\mc{Z}(x)$ is $\kappa_2(x)$ that is piecewise affine on a
polytopic partition $\{X_{21},X_{22},X_{23}\}$ of
$\mc{P}_2^\mc{X}=\mc{X}$. The sets $X_{ij}$ and the optimal control
laws $\kappa_1\bb$ and $\kappa_2\bb$ are shown in the Figure. At each
$x \in \mc{X}$, there are two candidates $\kappa_1(x)$ and
$\kappa_2(x)$ for the optimal control $\kappa(x)$ for the original
problem $\mathbb{P}_\textrm{min}(x)$. Theorem 1 resolves this
difficulty; at all $z$ on the boundary between $P_1^\mc{Z}$ and
$P_2^\mc{Z}$, $\mc{S}^\mc{Z}(z)=\{1,2\}$ since both polytopes are
active. Hence on the boundary, a control $u$ is optimal for the
original problem $\mathbb{P}_\text{min}(x)$ if and only if it is
optimal for both problems $\mc{P}_1^\mc{Z}(x)$ and
$\mc{P}_2^\mc{Z}(x)$.  At all $x \in X_{11}\cup X_{12}$ (except at its
intersection with $X_{13}$), only $(x,\kappa_2(x))$ lies on the
boundary between $P_1^\mc{Z}$ and $P_2^\mc{Z}$ ($(x,\kappa_1(x))$ does
not lie on this boundary). Thus, in $X_{11}\cup X_{12}$, the optimal
control is not $\kappa_2(x))$; it is $\kappa(x) = \kappa_1(x)$.  At
all $x \in$$X_{13}\cap X_{21}$ both $(x,\kappa_1(x))$ and
$(x,\kappa_2(x))$ lie on the boundary between $P_1^\mc{Z}$ and
$P_2^\mc{Z}$ so that the optimal control here is
$\kappa(x)=\kappa_1(x)=\kappa_2(x)$. Finally, at all $x \in$
$X_{22}\cup X_{23}$ (except at its intersection with $X_{21}$), only
$(x,\kappa_1(x))$ lies on the boundary between $P_1^\mc{Z}$ and
$P_2^\mc{Z}$ ($(x,\kappa_2(x))$ does not lie on this boundary); thus,
in $X_{22}\cup X_{23}$, the optimal control is $\kappa(x) =
\kappa_2(x)$. Hence $\kappa\bb$ is completely defined on $\mc{X}$.
This procedure avoids the overlapping of sets that results in previous
analysis of this problem.

\subsection{The maximization  subproblem $\mathbb{P}_\text{max}$}
In the minimization subproblem \eqref{eq:3.1}, the function $V\bb$
being minimized is convex in both $x$ and $u$. In contrast, in the
maximization subproblem \eqref{eq:3.2}, the function $V\bb$ being
maximized is convex in $z$ and (under suitable conditions) concave in
$w$.  Hence we proceed somewhat differently.

\begin{proposition}
  \label{p:5}
  Suppose $V:\Phi \tends \reals$ is such that $z \mapsto V(z,w)$ is
  strictly convex and continuous for each $w \in W$, $w \mapsto
  V(z,w)$ is strictly concave and continuous for each $z$ in $\mc{Z}
  \eqbyd \proj_Z(\Phi)$, and that $\Phi$ is a \emph{polytope} with a
  non-empty interior. Then, for all $(z,w) \in \Phi$, the solution
  $\nu(z)$ to $\mathbb{P}_\mathrm{max}(z)$ exists and is unique.
  Moreover, the value function $J^0\bb$ is strictly convex and
  continuous with domain $\mc{Z}$, and $\nu\bb$ is continuous in
  $\mc{Z}$.
\end{proposition}
\begin{proof}
  Since $J^0\bb$ is the maximum of a set of strictly convex and
  continuous functions, it is also strictly convex and continuous. The
  existence and uniqueness of $\nu(z)$ for each $z \in \mc{Z}$ follows
  from the strict concavity and  continuity of $w \mapsto V(z,w)$
  and the compactness of $W$. The continuity of $\nu\bb$ follows from
  the uniqueness of $\nu(z)$ at each $z$ (Theorem 5.4.3 in
  \cite{polak:1997}).
\end{proof}
To obtain a parametric solution to $\mathbb{P}_\mathrm{max}$, we
assume that $V\bb$ is continuous piecewise quadratic on a polytopic
partition $\mc{P}^\Phi =\{P_i^\Phi \mid i \in \mc{J}^\Phi\}$ of $\Phi
\eqbyd \mc{Z} \times W$ (in the absence of additional restrictions,
$\mc{Z} = X \times U$ so both $\mc{Z}$ and $\Phi$ are polytopic).
For each $(z,w) \in \Phi$, let $S^\Phi(z,w)$, the index set of active
polytopes at $(z,w)$, be defined by
\begin{equation}
  \label{eq:3.23}
  S^\Phi(z,w) \eqbyd \{i \in \mc{J}^\Phi \mid (z,w) \in P_i^\Phi\}
\end{equation}
so that $S^\Phi(z,w)$ is the set of indices of active polytopes at
$(z,w)$.  Similarly, for each $z \in \mc{Z} \eqbyd \proj_Z(\Phi)$, let
$S_\Phi^0(z)$ be defined by
\begin{equation}
  \label{eq:3.24}
  S_\Phi^0(z) \eqbyd S^\Phi(z,\nu(z))
\end{equation}
where
\begin{equation}
  \label{eq:3.25}
  \nu(z) \eqbyd \arg \max_{w} \{V(z,w) \mid w \in W\}
\end{equation}
so that $S_\Phi^0(z)$ is the index set of polytopes active at the
solution $\nu(z)$ of $\mathbb{P}_\text{max}(z)$. For each $i \in
\mc{J}^\Phi$, each $z \in \mc{Z}$, we define the simpler problem
$\mathbb{P}_i(z)$ defined by
\begin{alignat}{1}
  \label{eq:3.26}
 \mathbb{P}_i(z): \qquad J_i^0(z) &= \max_w \{V_i(z,w) \mid (z, w) \in
  P_i^\Phi\}\\
  \label{eq:3.27}
  \nu_i(z) &= \arg \max_w \{V_i(z,w) \mid  (z, w) \in  P_i^\Phi\}
\end{alignat}
where $V(z,w)=V_i(z,w)$ on $P_i^\Phi$ and $V_i\bb$ is quadratic.  For
each $i$, problem $\mathbb{P}_i^\Phi(z)$ is a quadratic program.
\begin{proposition}
  \label{p:6}
  Suppose $V: \Phi\tends \reals$ is strictly concave (hence
  continuous) in $w$ for each $z \in \proj_Z\Phi$, and
  that $\Phi$ is a polytope with a polytopic partition $\{P_i^\Phi
  \mid i \in \mc{J}^\Phi\}$ such that, for each $i \in \mc{J}^\Phi$,
  $P_i^\Phi$ has a non-empty interior. Then, for all $z \in
  P_i^\mc{Z}=\mathrm{Proj}_Z P_i^\Phi$, all $i \in \mc{J}^\Phi$, the
  solution $\nu_i(z)$ to $\mathbb{P}_i(z)$ exists and is unique.
  Moreover, the value function $J_i^0\bb$ is strictly convex (hence
  continuous) with domain $P_i^\mc{Z}$, and $\nu_i\bb$ is continuous at
  any $z \in P_i^\mc{Z}$.
\end{proposition}
The proof of Proposition \ref{p:6} is similar to the proof of
Proposition \ref{p:5}. The relation between the solution to
$\mathbb{P}_\mathrm{max}$ and the solutions to the subproblems
$\mathbb{P}_i^\Phi$, $i \in \mc{J}^\Phi$  is given in the next
result.
\begin{proposition}
  \label{p:7}
  Suppose $V: \Phi \tends \reals$ is continuous piecewise quadratic,
  strictly convex in $z$ and strictly concave in $w$ and is continuous
  piecewise quadratic in a polytopic partition $\mc{P}^\Phi
  =\{P_i^\Phi \mid i \in \mc{J}^\Phi\}$ of the polytope $\Phi$.  Then
  $w$ is optimal for the maximization problem
  $\mathbb{P}_\mathrm{max}(z)$ if and only if $w$ is optimal for the
  problems $\mathbb{P}_i(z)$ $(\nu(z)=\nu_i(z))$ for all $i \in
  {S}_{\Phi}^0(z)$.
\end{proposition}
The proof of Proposition \ref{p:7} is similar to the proof of
Proposition \ref{p:3}. We now exploit the continuous piecewise quadratic nature
of $V\bb$. For each $i \in \mc{J}^\Phi$, subproblem $\mathbb{P}_i(z)$
may be expressed as:
\begin{equation}
  \label{eq:3.28}
  \mathbb{P}_i(z): \qquad J_i^0(z) = \max_w \{V_i(z,w) \mid M_iw
 \le N_iz +p_i\}
\end{equation}
($V_i\bb$ is quadratic) for some $M_i,N_i,p_i$, each matrix (vector)
having $r_i$ rows.  If we assume that the constraints indexed by $I
\subseteq \mc{I}_{r_i}$ are active, then $ \mathbb{P}_i(z)$ is replaced
by the simpler, equality constrained, problem
\begin{equation}
  \label{eq:3.29}
  \mathbb{P}_{i,I}(z): \qquad J_i^0(z) = \max_w \{V_i(z,w) \mid
  M^j_iw = N^j_iz+p_i^j,\ j \in I \}
\end{equation}
where the superscript $j$ on matrix (or vector) denotes the
$j^\text{th}$ row of the matrix (or vector). The solution to this
problem is
\begin{alignat}{1}
  \label{eq:3.30}
  J_{i,I}^0(z) &=(1/2)z'Q_{i,I}z+q_{i,I}'z+s_{i,I}\\
  \label{eq:3.31}
  \nu_{i,I}(z) &= K_{i,I}z+k_{i,I}
\end{alignat}
Let $M_i^I$ denote the matrix the rows of which are $M_i^j$, $j \in
I$.  Let $PC_{i,I}\eqbyd \{(M_i^I)'\lambda \mid \lambda \ge 0\}$
$\subseteq \reals^p$ denote the polar cone at $0$ to the cone
$\mc{F}_{i,I} = \{h \mid M^j_i h \le 0,\ j \in I\}$ of feasible
directions $h$ for problem $\mathbb{P}_i(z)$ at $w=\nu_{i,I}(z)$; the
polar cone depends \emph{solely} on $I$, the set of active
constraints; it does not depend on the parameter $z$.  The following
result holds \cite{mayne:rakovic:2003a}:
\begin{proposition}
  \label{p:8}
  The affine control law $\nu_{i,I}\bb$ is optimal for problem
  $\mc{P}_i(z)$ at all $z$ in the polytope $Z_{i,I}$
  defined by
\begin{equation}
  \label{eq:3.32}
Z_{i,I} \eqbyd \left\{ z \in \mc{Z}\ \pmb{\Bigg\vert} \begin{array}{rcl}
  M_i^j(K_{i,I}z+k_{i,I}) &\le& N_i^jz+p_i^j,\ j \in
  \mc{I}_{r_i} \setminus I\\
   \nabla_w V_i(z,\nu_{i,I}(z))&\in&PC_{i,I} \end{array} \right\}
\end{equation}
\end{proposition}
(the restriction $z \in \mc{Z}=X \times U$ is included in the
definition of $Z_{i,I}$ since it is not included as a constraint in
$\mathbb{P}_i(z)$).  As before, we have to `marry' a set of polytopes
$Z_{i,I}$ for all $i$ such that polytope $P_i^\Phi$ is active.
Suppose, for each $i \in \mc{J}^\Phi$, polytope $P_i^\Phi$ is defined
by
\begin{equation}
  \label{eq:3.33}
  P_i^\Phi \eqbyd \{(z,w) \mid M_iw \le N_iz+p_i\}
\end{equation}
for some $M_i,N_i,p_i$ each having $r_i$ rows.  For each $i\in\mc{J}^\Phi$,
let $P_i^\mc{Z}$ denote the polytope defined by
\begin{equation}
  \label{eq:3.34}
  P_i^\mc{Z} \eqbyd \{z  \mid \exists\ w \text{ s.t. } (z,w) \in P_i^\Phi\}
  = \proj_Z(P_i^\Phi)
\end{equation}
For each $z \in P_i^\mc{Z}$, each $i \in \mc{J}^\Phi$, the set of
active constraints for $\mathbb{P}_i(z)$ is
\begin{equation}
  \label{eq:3.35}
  I_i^0(z) \eqbyd \{ j \in \mc{I}_{r_i} \mid M_i^j\nu_i(z) = N_i^jz+p_i^j\}
\end{equation}
where $M_i^j$ is the $j^\text{th}$ row of $M_i$, $N_i^j$ the
$j^\text{th}$ row of $N_i$, and $p_i^j$ the $j^\text{th}$ row of
$p_i$.  Also, for each $\bar{z} \in \mc{Z}$,  let the polytope
$Z_i(\bar{z}) \subseteq P_i^\mc{Z}$ be defined by
\begin{equation}
  \label{eq:3.36}
  Z_i(\bar{z})\eqbyd   Z_{i,I_i^0(\bar{z})}
\end{equation}
where, for each index set $I\subseteq \mc{J}^\Phi$, $Z_{i,I}$ is
defined by \eqref{eq:3.32}.  The polytope $Z(\bar{z})$ that figures in
the parametric solution of $\mathbb{P}_\text{max}(z)$ is defined, for
each $\bar{z} \in \mc{Z}$, by
\begin{equation}
  \label{eq:3.37}
  Z(\bar{z}) \eqbyd \cap \{Z_i(\bar{z}) \mid i \in
  S_{\Phi}^0(\bar{z})\}
\end{equation}
Clearly $\bar{z} \in Z(\bar{z})$ and $Z(z')=Z(\bar{z})$ for all $z'
\in Z(\bar{z})$. For each $\bar{z} \in \mc{Z}$, let the functions
$J_{\bar{z}}^0(\cdot)$ and $\nu_{\bar{z}}(\cdot)$ be defined on
$Z(\bar{z})$ by
\begin{alignat}{1}
  \label{eq:3.38}
  J_{\bar{z}}^0(z) &\eqbyd J_i^0(z),\ \forall i \in
  S_{\Phi}^0(\bar{z})\\
  \label{eq:3.39}
  \nu_{\bar{z}}(z) &\eqbyd \nu_i(z),\ \forall i \in
  S_{\Phi}^0(\bar{z})
\end{alignat}
The domain of each function is $Z(\bar{z})$; that the functions are
well defined follows from Proposition \ref{p:7}, \eqref{eq:3.35} and
\eqref{eq:3.36} which show that $\nu_i(z)=\nu_j(z)$ for all $i,j \in
S_\Phi^0(\bar{z})$, all $z \in Z(\bar{z})$.  Summarizing, we have:
\begin{theorem}
  \label{th:2}
  Suppose $V:{\Phi} \tends \reals$ is continuous piecewise quadratic
  and strictly convex on a polytopic partition
  $\mc{P}^\Phi=\{P_i^\Phi \mid i \in \mc{J}^\Phi\}$ of the polytope
  ${\Phi}$ and that, for each $i \in \mc{J}^\Phi$, $P_i^\Phi$ has a
  non-empty interior. Then the value function $J^0\bb$ is continuous
  piecewise quadratic and strictly convex on a polytopic partition
  $\mc{P}^\mc{Z}=\{P_i^\mc{Z} \mid i \in \mc{J}^\mc{Z}\}$ of $\mc{Z}
  \eqbyd \{z \mid \exists w \in W \mathrm{\ s.t.\ }(z,w) \in \Phi\}$
  $=\proj_Z(\Phi)$. The maximizer $\nu\bb$ is piecewise affine on
  $\mc{P}^\mc{Z}$. The polytopes $P_i^\mc{Z}$ are each of the form
  $Z(\bar{z})$ for some $\bar{z} \in \mc{Z}$; the value function and
  optimal control law satisfy $J^0(z) = J_{\bar{z}}^0(z)$ and
  $\nu(z)=\nu_{\bar{z}}(z)$ for all $z \in Z(\bar{z})$ and some
  $\bar{z} \in \mc{Z}$.
\end{theorem}
The proof of this result follows from Propositions \ref{p:3} and
\ref{p:4} and the discussion above.  As stated above, in the absence
of further restrictions, $\Phi = X \times U \times W$. However, in our
use of this result, the cost function $V\bb$ has the form
$V(z,w)=\ell(z,w) +V^0(Fz+Gw)$ where $F \eqbyd [A,B]$ and $V^0(x)$ may
be known only on a compact subset $\mc{X}$ of $\reals^n$; in this case
$\Phi = \{(z,w) \in X \times U \times W \mid Fz+Gw \in \mc{X}\}$. If
$A$ is invertible (which we assume) or $X$ is compact, then $\Phi$ is
a (compact) polytope with a polytopic partition.

\section{$H_\infty$ control; no state constraints}
In this section, we consider the $H_\infty$ constrained optimal control
problem when the only constraints are $u \in U$ and $w \in W$,  i.e.
$X=X_f=\reals^n$.  In this case, the dynamic programming equations
\eqref{eq:2.1} - \eqref{eq:2.4} simplify and are replaced by the
conventional dynamic programming equations:
\begin{alignat}{1}
  \label{eq:4.1}
  V_j^0(x) &= \min_{u \in U} \max_{w \in W}\{\ell(x,u,w)+
           V_{j-1}^0(f(x,u,w))\}\\
  \label{eq:4.2}
  \kappa_j(x) &= \arg\min_{u \in U}  \max_{w \in
  W}\{\ell(x,u,w)+V_{j-1}^0(f(x,u,w))\}
\end{alignat}
with boundary condition
\begin{equation}
  \label{eq:4.3}
  V_0^0(x) = V_f(x)
\end{equation}
The domain $X_j$ of $V_j^0$ now satisfies $X_j=\reals^n$, for all $j
\ge 0$ so the recursion equation \eqref{eq:2.3} for $X_j$ is not
required. The recursion equations may be rewritten in the form
\begin{alignat}{1}
  \label{eq:4.4}
  V_j^0(x) &= \min_{u \in U} J_{j-1}^0(x,u) \\
  \label{eq:4.5}
  J_{j-1}^0(z) &\eqbyd  \max_{w \in W}
  \{\ell(z,w)+V_{j-1}^0(f(z,w))\}\\
  \label{eq:4.6}
  \kappa_j(x) &= \arg\min_{u \in U} J_{j-1}^0(x,u)\\
  \label{eq:4.7}
  \nu_j(z) &= \arg\max_{w \in W}
  \{\ell(z,w)+V_{j-1}^0(f(z,w))\}
\end{alignat}
where $z=(x,u)$.
The important feature of this problem is that the constraint $u\in U$
in subproblem \eqref{eq:4.4} has the same simple form as the
constraint $w \in W$; this permits us to obtain stronger properties
for the value functions $V_j^0\bb$ and $J_j^0\bb$, $j \ge 0$. The
prototype problems for \eqref{eq:4.4} and \eqref{eq:4.5} are, respectively:
\begin{alignat}{1}
  \label{eq:4.8}
  \mathbb{P}_\mathrm{min}(x): \qquad V^0(x) &= \min_{u} \{J(x,u) \mid
  u \in U\}\\
  \label{eq:4.9}
  \mathbb{P}_\mathrm{max}(z): \qquad J^0(z) &= \max_{w} \{V(z,w)  \mid
  w \in W\}
\end{alignat}
in which $V^0\bb$ and $J\bb$ replace, respectively, $V_j^0\bb$ and
$J_{j-1}^0\bb$ in \eqref{eq:4.4} and $J^0\bb$ and $V\bb$ replace,
respectively, $J_{j-1}^0\bb$ and $\ell\bb+V_{j-1}(f(\cdot))$ in
\eqref{eq:4.5}.  Since $J\bb$ is convex in $x$ and $V\bb$ is (under
appropriate conditions) concave in $w$, their respective value
functions have identical properties ($\min_u \{J(x,u) \mid u \in U\} =
-\max_u \{-J(x,u) \mid u \in U\}$).

\begin{proposition}
  \label{p:9}
  Suppose that $J\bb$ in $\mathbb{P}_\mathrm{min}$ $(V\bb$ in
  $\mathbb{P}_\mathrm{max})$ is continuously differentiable, strictly
  convex in $u$ (strictly concave in $w$) and that $U$ ($W$) is
  compact. Then the value function $V^0\bb$ of
  $\mathbb{P}_\mathrm{min}$ $(J^0\bb$ of $\mathbb{P}_\mathrm{max})$
  is continuously differentiable.
\end{proposition}
\begin{proof}
  It is only necessary to consider $\mathbb{P}_\mathrm{max}$ and
  establish the continuous differentiability of $J^0\bb$.  Since $W$,
  being constant, is continuous in $z$, the continuity of the value
  function $J^0\bb$ follows from the maximum theorem (e.g.  Theorem
  5.4.1 in \cite{polak:1997}).  Since the function $w \mapsto V(z,w)$
  is strictly concave for all $z$, the maximizer $\nu(z)$ is unique (a
  singleton) for each $z$; by the same maximum theorem, $\nu\bb$ is
  continuous.  Since $V\bb$ is continuously differentiable and $W$ is
  compact, and the maximizer $\nu(z)$ is unique and continuous, it
  follows from the proof of Theorem 5.4.7 in \cite{polak:1997} that
  the directional derivative of $J^0\bb$ satisfies
  \begin{equation}
    \label{eq:4.10}
     d J^0(z;h) = (\partial/\partial z)V(z,\nu(z)) h
  \end{equation}
  at any $z$, any direction $h$.  Hence $J^0\bb$ is Gateau
  differentiable at any $z \in \mc{Z}$ with Gateau derivative
  $G(z)=(\partial/\partial z) V(z,\nu(z))$. Since $G\bb$ is
  continuous, $J^0\bb$ is continuously (Frechet) differentiable in
  $\mc{Z}$ with derivative $(\partial/\partial z)J^0(z)=
  (\partial/\partial z) J(z,\nu(z))$ \cite{vainberg:1974}.
\end{proof}

Although the domain of the value functions $V_j^0\bb$ is $\reals^n$
and that of the value functions $J_j^0\bb$ is $\reals^{n+m}$, we
restrict attention in this section to polytopic subsets of these
domains. With this caveat, we now show that there exists a $\gamma>0$
such that $V\bb$ in $\mathbb{P}_\text{max}$, which represents
$V_{j-1}^0\bb$ in \eqref{eq:4.5} and therefore has the form
$V(z,w)=\ell(z,w)+V^0(Fz+Gw)$ where $F\eqbyd[A,B]$, is strictly
concave in $w$.
\begin{proposition}
  \label{p:10}
  Let $\mc{X}$ be a polytope in $\reals^n$ containing the origin in
  its interior.  Suppose $V\bb$ is defined by $V(z,w) \eqbyd
  \ell(z,w)+V^0(Fz+Gw)$ where $V^0\bb$ is continuously differentiable
  and continuous piecewise quadratic on a polyhedral partition
  $\mc{P}^\mc{X}=\{P_i^\mc{X} \mid i \in \mc{J}^\mc{X}\}$ of $\mc{X}$.
  Then $V\bb$ is continuously differentiable and continuous piecewise
  quadratic on a polyhedral partition $\mc{P}^\Phi=\{P_i^\Phi \mid i
  \in \mc{J}^\Phi\}$ of the polyhedron $\Phi \eqbyd\{(z,w) \in \reals^n
  \times U \times W \mid Fz+Gw \in \mc{X}\}$ and there exists a
  $\gamma^\ast>0$ such that $V\bb$ is strictly concave in $w$ for each
  $z$ in $\mc{Z}\eqbyd\proj_Z\Phi$ and all $\gamma\ge \gamma^\ast$.
\end{proposition}
\begin{proof}
  The continuous differentiability of $V\bb$ follows from the
  continuous differentiability of $\ell\bb$ and $V^0\bb$.
  Take any two points $w_1,w_2$ in $W$. For all $\lambda \in [0,1]$,
  let $w_\lambda \eqbyd w_1+\lambda(w_2-w_1)$, and, for each $z \in
  \mc{Z}$, let the real valued function $\phi\bb$ be defined on
  $[0,1]$ by $\phi(\lambda) \eqbyd V(z,w_\lambda)$.  Suppose that
  $V^0(x) = (1/2)x'Q_ix+q_i'x+r_i$ in $P_i^\mc{X}$ (for each $i \in
  \mc{J}^\mc{X}$). Then
  \begin{alignat*}{1}
  V(z,w) &= (1/2)(Fz+Gw)'Q_i(Fz+Gw) +q_i'(Fz+Gw)+r_i +\ell(z,w)\\
         &= -(1/2)w'(\gamma^2I-G'Q_iG)w+ b_i'w+c_i
  \end{alignat*}
  on the polyhedron $P_i^\Phi = \{(z,w)\in \reals^n \times U \times W
  \mid Fz+Gw \in P_i^\mc{X}\}$, where $b_i$ and $c_i$ depend on $z$.
  For any $\varepsilon > 0$, there exists a $\gamma^\ast>0$ such that
  $C_i \eqbyd \gamma^2I-G'Q_iG\ge \varepsilon I$ for all $\gamma \ge
  \gamma^\ast$, all $i \in \mc{J}^\mc{X}$. The function $\phi\bb$ is
  continuously differentiable and satisfies:
  \begin{alignat*}{1}
   \phi(\lambda) &= -(1/2)(h'C_ih) \lambda^2+b_i\lambda+c_i\\
   \phi'(\lambda) &= -(h'C_ih) \lambda+b_i
  \end{alignat*}
  for all $\lambda\in[0,1]$ such that $Fz+Gw_\lambda \in P_i^\Phi$.
  Since $\phi'\bb$ is continuous, $\phi'\bb$ is strictly decreasing if
  $\gamma\ge \gamma^\ast$ . It follows, by a trivial modification to
  the proof of Theorem 4.4 in \cite{rockafellar:1970}, that
  $\phi(\lambda) > \phi(0)+\lambda(\phi(1) -\phi(0))$ for all $\lambda
  \in (0,1)$ which establishes the strict concavity of $\phi\bb$ and,
  hence, of $w \mapsto V(z,w)$ if $\gamma\ge \gamma^\ast$. That $V\bb$
  is piecewise quadratic on a polyhedral partition of $\mc{P}^\Phi$
  follows, with minor amendments, from the proofs of Proposition
  \ref{p:8} and Theorem \ref{th:2}.
\end{proof}
We can now establish the main result of this section, characterization
of the solution to the constrained $H_\infty$ problem when
$X=X_f=\reals^n$.  We characterize the value functions $V_j^0$ on polytopic
subsets $X_j$ of the true domain $\reals^n$ by assuming that the
terminal cost function $V_f\bb$ is known only in a polytopic subset
$X_0$ of $\reals^n$.

\begin{theorem}
  \label{th:3}
  Suppose $V_f\bb$ is continuously differentiable, strictly convex,
  and continuous piecewise quadratic on a polytopic partition
  $\mc{P}_0^X$ of a polytope $X_0 \subset \reals^n$. Then, there
  exists a $\gamma>0$ such that, for each $j \ge 0$, there exists a
  polyhedron ${X}_j$ on which the value function $V_j^0\bb$ is
  continuously differentiable, strictly convex, and continuous
  piecewise quadratic on a polyhedral partition $\mc{P}_j^X$ of
  ${X}_j$, and the optimal control law $\kappa_j\bb$ is continuous and
  piecewise affine on the same polyhedral partition $\mc{P}_j^X$ of
  ${X}_j$.
\end{theorem}
\begin{proof}
  Suppose $V_{j-1}^0\bb$ is continuously differentiable, strictly
  convex, and continuous piecewise quadratic on a polyhedral partition
  $\mc{P}^{X}_{j-1}$ of a polyhedron ${X}_{j-1}$ if $\gamma\ge
  \gamma_{j-1}$. Then, by Proposition \ref{p:10}, there exists a
  $\gamma_j \ge \gamma_{j-1}$ such that $(z,w) \mapsto
  \ell(z,w)+V_{j-1}^0(f(z,w))$ is strictly concave in $w$,
  continuously differentiable and continuous piecewise quadratic on a
  polyhedral partition $\mc{P}^\Phi_{j-1}$ of a polyhedron $\Phi_{j-1}
  =\{(z,w) \in \reals^n \times U \times W \mid Fz+Gw \in {X}_{j-1}\}$.
  By Proposition \ref{p:9}, the value function $J_{j-1}^0\bb$ is then
  continuously differentiable and, by Theorem \ref{th:2},
  $J_{j-1}^0\bb$ is continuous piecewise quadratic and strictly convex
  on a polyhedral partition $\mc{P}_{j-1}^{Z}$ of a polyhedron
  ${Z}_{j-1}=\proj_Z \Phi_{j-1}$ (and the disturbance law
  $\nu_{j-1}\bb$ is continuous and piecewise affine on the same
  polytopic partition).  Then, by Proposition \ref{p:9}, $V_j^0\bb$ is
  continuously differentiable and, by Theorem \ref{p:1}, $V_j^0\bb$ is
  strictly convex and continuous piecewise quadratic on a polyhedral
  partition $\mc{P}_j^{X}$ of a polyhedron ${X}_j$ (and the optimal
  control law $\nu_{j-1}\bb$ is continuous and piecewise affine on the
  same polyhedral partition).
 \end{proof}

\section{$H_\infty$ control; state and control constraints}

In this section, we consider the $H_\infty$ constrained optimal
control problem when the constraints are $u \in U$, $w \in W$ , $x\in
X$ and the terminal constraint $x_N \in X_f$. The dynamic programming
solution of the $H_\infty$ problem requires the repeated solution of
the two prototype problems $\mathbb{P}_\text{min}(x)$ and
$\mathbb{P}_\text{max}(z)$ defined in \eqref{eq:3.1} and
\eqref{eq:3.2} which we rewrite in the form:

\begin{alignat}{1}
  \label{eq:5.1}
  \mathbb{P}_\mathrm{min}(x): \qquad V^0(x) &= \min_{u} \{J(x,u)
  \mid u \in \mc{U}(x)\}\\
  \label{eq:5.2}
  \mathbb{P}_\mathrm{max}(z): \qquad J^0(z) &= \max_{w} \{V(z,w)  \mid
  w \in W\}
\end{alignat}
where the set-valued function $\mc{U}\bb$ is defined, for all $x \in
\mc{X} \eqbyd \proj_X(Z)$, by
\begin{equation}
  \label{eq:5.3}
  \mc{U}(x) \eqbyd \{u \mid (x,u) \in Z\}
\end{equation}
The presence of state constraints complicates the solution of the
$H_\infty$ problem considerably.  The extra complexity arises in the
solution of $\mathbb{P}_\text{min}(x)$ since the control constraint $u
\in \mc{U}(x)$ is now dependent on the parameter $x$ in contrast to
the simple constraint $u \in U$ when no state constraints are present.
The dependency of the constraint on $x$ can cause the gradient of the
value function $V^0\bb$ in problem $\mathbb{P}_\text{min}(x)$ to be
discontinuous even if the function $J\bb$ being minimized is
continuously differentiable; Proposition 9 is no longer necessarily
true for problem $\mathbb{P}_\text{min}(x)$.  However, there do exist
conditions under which this result is true.

\subsection{Particular case}

Assume that $J\bb$ in \eqref{eq:3.1} is continuously differentiable
and continuous piecewise quadratic on a polytopic partition
$\mc{P}^\mc{Z}=\{P_i^\mc{Z}, i \in \mc{J}^\mc{Z}\}$ of $\mc{Z}$.  We
show below, despite the fact that the constraint set $\mc{U}(x)$ now
depends on the parameter $x$, that the value function $V^0\bb$ for
$\mathbb{P}_\text{min}$ is, under certain further assumptions,
continuously differentiable in $\mc{X} \eqbyd \proj_X(\mc{Z})$.  To do
this, we first consider, as in \S3.1, the simpler problems
$\mathbb{P}_i(z)$, $i \in \mc{J}^\mc{Z}$, defined by \eqref{eq:3.8}.
For each $i$, problem $\mathbb{P}_i(z)$ is a quadratic program, with a
value function $V_i^0\bb$ that is continuous piecewise quadratic on a polytopic
partition of the polytope $P_i^\mc{X} \eqbyd \proj_X(P_i^\mc{Z})$, and
may be written in the form
\begin{equation}
  \label{eq:5.4}
  \mathbb{P}_i(x): \qquad V_i^0(x) = \min_u \{J_i(x,u) \mid
  x \in  \mc{U}_i(x)\}
\end{equation}
where
\begin{equation}
  \label{eq:5.5}
  \mc{U}_i(x) \eqbyd \{u \mid (x,u) \in P_i^\mc{Z}\}=
       \{u \mid M_iu \le N_ix+p_i\}
\end{equation}
It is known (see Proposition \ref{p:4} and Theorem \ref{th:1}) that
the value function $V_i^0\bb$ is continuous piecewise quadratic, being
quadratic on polytopes $X_{i,I}$, each polytope characterized by a
set $I \subseteq \mc{I}_{r_i}$ of active constraints, where $r_i$ is
the number of rows of $M_i$; the sets $X_{i,I}$, $I \subseteq
\mc{I}_{r_i}$ (excluding sets with no interior) constitute a polytopic
partition of the polytope $P_i^\mc{X}$. We require the following
result which is proved in the appendix.
\begin{proposition}
  \label{p:11}
Suppose (i), $J_i\bb$ is continuously differentiable, (ii)
$P_i^\mc{Z}$ has an interior, and, (iii) for any two adjacent
polytopes,  $X_{i,I_1}$ and $X_{i,I_2}$ say, in
the polytopic partition of $P_i^\mc{X}$, either $I_1 \subseteq I_2$ or
$I_2 \subseteq I_1$.  Then,  $V_i^0\bb$ is continuously differentiable
in  $P_i^\mc{X}$.
\end{proposition}
We establish next the continuous differentiability of
the value function $V^0\bb$ for $\mathbb{P}_\text{min}$.
\begin{theorem}
  \label{th:4}
Suppose that $J\bb$ in \eqref{eq:3.1} is continuously differentiable
and continuous piecewise quadratic on a polytopic partition
$\mc{P}^\mc{Z}=\{P_i^\mc{Z}, i \in \mc{J}^\mc{Z}\}$ of $\mc{Z}$ and
that hypotheses (ii) and (iii) of Proposition \ref{p:11} are
satisfied (hypothesis (i) is satisfied automatically)
for each problem $\mathbb{P}_i$, $i \in \mc{J}^\mc{Z}$.  Then
$V^0\bb$, the value function for $\mathbb{P}_\text{min}$, is
continuously differentiable in $\mc{X}$.
\end{theorem}
\begin{proof}
  It follows from Theorem \ref{th:1}, that, for each $i \in
  \mc{J}^\mc{Z}$, $V^0(x)=V_i^0(x)$ for all $x$ in the interior of
  each polytope $X_{i,I}$ in the polytopic partition of $P_i^\mc{X}$.
  It follows from Proposition \ref{p:11}, that $V^0\bb$ is
  continuously differentiable in $P_i^\mc{X}$ for each $i \in
  \mc{J}^\mc{Z}$. Consider next a point $\bar{x}$ on the boundary
  between two polytopes $P_i^\mc{X}$ and $P_j^\mc{X}$; clearly $i$ and
  $j$ both lie in $S_\mc{Z}^0(\bar{x})$ and, from Theorem \ref{th:1},
  $V^0(x) = V_i^0(x) =V_j^0(x)$ for all $x \in X(\bar{x})$, so that
  $V^0\bb$ is continuously differentiable in $X(\bar{x})$ and, hence,
  on all boundaries between polytopes in the polytopic partition of
  $\mc{X}$.
\end{proof}

\begin{theorem}
  \label{th:5}
  Suppose $V_f\bb$ is continuously differentiable, strictly convex,
  and continuous piecewise quadratic on a polytopic partition
  $\mc{P}_0^X$ of a polytope $X_0 \subset \reals^n$. Then, there
  exists a $\gamma>0$ such that, for each $j \ge 0$, there exists a
  polyhedron ${X}_j$ on which the value function $V_j^0\bb$ is
  continuously differentiable, strictly convex, and continuous
  piecewise quadratic on a polytopic partition $\mc{P}_j^X$ of
  ${X}_j$, and the optimal control law $\kappa_j\bb$ is continuous and
  piecewise affine on the same polytopic partition $\mc{P}_j^X$ of
  ${X}_j$.
\end{theorem}
\begin{proof}
  Suppose $V_{j-1}^0\bb$ is continuously differentiable, strictly
  convex, and continuous piecewise quadratic on a polytopic partition
  $\mc{P}^{X}_{j-1}$ of a polytope ${X}_{j-1}$ if $\gamma\ge
  \gamma_{j-1}$. Then, by Proposition \ref{p:10}, there exists a
  $\gamma_j \ge \gamma_{j-1}$ such that $(z,w) \mapsto
  \ell(z,w)+V_{j-1}^0(f(z,w))$ is strictly concave in $w$,
  continuously differentiable and continuous piecewise quadratic on a
  polytopic partition $\mc{P}^\Phi_{j-1}$ of a polytope $\Phi_{j-1}
  =\{(z,w) \in X \times U \times W \mid Fz+Gw \in {X}_{j-1}\}$.
  By Proposition \ref{p:9}, the value function $J_{j-1}^0\bb$ is then
  continuously differentiable and, by Theorem \ref{th:2},
  $J_{j-1}^0\bb$ is continuous piecewise quadratic and strictly convex
  on a polytopic partition $\mc{P}_{j-1}^{Z}$ of a polytope
  ${Z}_{j-1}=\proj_Z \Phi_{j-1}$ (and the disturbance law
  $\nu_{j-1}\bb$ is continuous and piecewise affine on the same
  polytopic partition).  Then, by Proposition \ref{p:9}, $V_j^0\bb$ is
  continuously differentiable and, by Theorem \ref{p:1}, $V_j^0\bb$ is
  strictly convex and continuous piecewise quadratic on a polytopic
  partition $\mc{P}_j^{X}$ of a polytope ${X}_j$ (and the optimal
  control law $\nu_{j-1}\bb$ is continuous and piecewise affine on the
  same polytopic partition).
 \end{proof}

\subsection{General case}

A simple characterization for the solution of the $H_\infty$ problem
with control and state constraints does not appear possible when the
simplifying assumption of \S5.1 is not made. Without this assumption,
the value function $V^0\bb$ for the minimization problem
$\mathbb{P}_\text{min}$ is not necessarily continuously differentiable
at the boundary between polytopes in the polytopic partition of
$\mc{X}$. Consequently, the objective function $V\bb$ (which has the
form $V=\ell+V^0$) in the maximization problem is not necessarily
concave, no matter how large $\gamma$ is chosen. The resultant cost
function $J\bb$ in the minimization problem is then piecewise
max-quadratic, i.e. it is continuous and equal to the maximum of a
finite number of quadratics in each polytope in a polytopic
partition of its domain.  It does not appear possible to obtain a
simple characterization for a problem with this structure.

\subsection{Illustrative example}

The partial value functions and optimal control laws can be computed
by solving the max and min subproblems associated (for each
subproblem) with each set of potentially active constraints and
each set of potentially active polytopes. Since the number of
these sets is combinatorial, a better procedure, employed in our
computations, is to select a state-control pair $z=(x,u)$ in the max
subproblem, determine the active constraint set $J$ or the active
index set of polytopes $s$, and then compute the corresponding
disturbance law $\nu\bb$ and the region $Z$ in which this control law
is optimal, using Theorem 2.  The procedure is then repeated for a new
value of $z$ not lying in the union of the sets $Z$ already computed.
Once the max subproblem computations are complete, a similar procedure
is applied to the min subproblem using Theorem 1.

Our numerical example is optimal min-max control of a constrained
second order system defined by:
\begin{equation}
\label{eq:7.1}
        x^+= \begin{array}{c}
                 \left[ \begin{array}{cc}
         1 & 1\\
         0 & 1\\
         \end{array}\right]x
         +\left[ \begin{array}{c}
         0.5\\
         1\\
         \end{array}\right]u
         + w\\
         \end{array}
\end{equation}
The state constraints are $x\in X:=\{x \mid  |x|_\infty\le10\}$.
The control constraint is $u\in U:=\{u \mid |u|\le1\}$.
The disturbance is bounded: $w\in W:=\{w \mid |w|_\infty\le 0.1\}$
The path cost function is quadratic with $Q=10I$, $R=1$ and $\gamma=100$.
The terminal cost $V_f(x)$ is quadratic $(1/2)x'P_fx$ with
\[
P_f=\left[ \begin{array}{cc}
         20.6143  &  5.9244\\
         5.9244   & 14.2329\\
         \end{array}\right]
\]
The terminal constraint set $X_f$ is defined  by the 4 inequalities:
$(-0.9849\   -0.3155)x\le 2.1526$; $(0.9489\    0.3155)x\le  2.1526$;
$(0.4369\     0.8995)x\le 0.7079$ and $(-0.4369\   -0.8995)x\le
0.7079$.

\begin{figure}[!ht]
\subfigure[Regions of x space, $N=1$]
{\includegraphics[width=5cm]{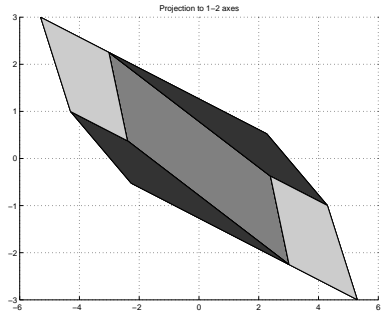}} \hspace{1in}
\subfigure[Regions of x space, $N=2$]
{\includegraphics[width=5cm]{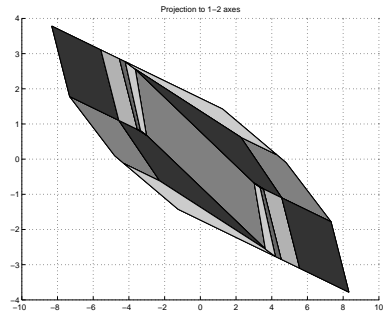}} \caption{2D Example}
\label{f1}
\end{figure}

The polytopic regions for $N=1$ and $N=2$ are shown in Figure
\ref{f1}.  At $N=2$, the state space is partitioned into 17 polytopes
in each of which the optimal control law is piecewise affine; the
state control space $\mc{Z}$ is partitioned into 5 polytopes in each
of which the optimal disturbance law is piecewise affine.

\section{$H_\infty$ receding horizon control}

\subsection{Introduction}

Since we make use, in this section, of the solution for infinite
horizon, linear unconstrained $H_\infty$ problem, we assume, in the
sequel, that $(A,B)$ is stabilizable and that $(C,A,B)$ has no zeros
on the unit circle where $Q = C'C$. Since $Q$ is assumed to be
positive definite, $(C,A)$ is detectable. These conditions, and the
fact that $R$ is assumed positive definite, ensure that the conditions
assumed in \cite{green:limebeer:1995}, Appendix B, are satisfied for
the full information case. Hence there exists a $\tilde{\gamma}>0$
such that a positive definite solution $P_f$ to the associated
(generalized) $H_\infty$ algebraic Riccati equation for all $\gamma >
\tilde{\gamma}$ and associated optimal control and disturbance laws
$u=K_ux$ and $w=K_w x$ respectively.  It is shown in
\cite{green:limebeer:1995} that, under these assumptions, the state
matrices $A_f \eqbyd A+BK_u$ and $A_c \eqbyd A+B K_u+GK_w$ are both
stable.

The terminal cost function $V_f\bb$ for the constrained $H_\infty$
control problem is defined by
\begin{equation}
  \label{eq:6.1}
  V_f(x) =(1/2)|x|_{P_f}^2.
\end{equation}
and satisfies
\begin{equation}
  \label{eq:6.2}
 V_f(A_cx)-V_f(x)+ \ell(x,K_ux,K_wx)  = 0.
\end{equation}
The terminal constraint set $X_f$ is chosen to be a disturbance
invariant set (if it exists) for the system $x^+= A_fx+Gw$, $A_f\eqbyd
A+BK_u$.  Any disturbance invariant set $X_f$ satisfies
\begin{equation}
  \label{eq:6.3}
   f(x,K_ux,W) \subseteq X_f\ \forall\ x \in X_f
\end{equation}
We assume that the set $W$ is sufficiently small, and that $\gamma$ is
sufficiently large, to ensure the existence of a disturbance
invariant set $X_f$ which satisfies
\begin{equation}
  \label{eq:6.4}
   X_f \subseteq X, \quad K_uX_f \subseteq U, \quad K_wX_f \subseteq
   W.
\end{equation}
That the last condition in \eqref{eq:6.4} can be satisfied follows
from the fact \cite{green:limebeer:1995} that $K_w \tends 0$ as
$\gamma \tends \infty$.  A suitable $X_f$ may be computed as follows:
if $W=\{w \mid C_ww \le c_w\}$, choose $\gamma$ such that $\tilde{X}
=\{x \mid C_w K_w x \le c_w\}$ is reasonably large; clearly $K_w
\tilde{X} \subseteq W$. Next, choose $X_f$ to be a disturbance
invariant set for $x^+= A_fx+Gw$ satisfying $X_f \subseteq X \cap
\tilde{X}$ and $K_u X_f \subseteq U$.

Since there does not exist a disturbance $w$ that can steer the system
outside $X_f$ given an initial state in $X_f$, the optimal policy for
$w$ in $X_f$ is $w=K_wx$. The closed loop system $x^+=A_\mathrm{c}x$,
$A_\mathrm{c} \eqbyd A+BK_u+GK_w$, is exponentially stable and the
controller $u=K_ux$ maintains the state in $X_f$ if the initial state
is in $X_f$. We observe that the solution of the infinite horizon
constrained $H_\infty$ problem (defined by \eqref{eq:2.9} with
$N=\infty$) satisfies:
\begin{equation}
   \label{eq:6.5}
   V_\infty^0(x)=V_f(x), \quad \kappa_\infty(x) = K_ux,\ \forall x \in
   X_f
\end{equation}
since, by \eqref{eq:6.4}, the control constraints are satisfied
everywhere in $X_f$ so that the solutions of the constrained and
unconstrained problems coincide.

\subsection{$H_\infty$ control: control constraints}

Since this problem, as stated in \S4, has no terminal constraint and
since $V_f\bb$ defined above is a local rather than a global Control
Lyapunov Function ($V_f\bb$ is valid in $X_f$), standard stability
results
\cite{nicolao:magni:scattolini:1998,mayne:rawlings:rao:scokaert:2000}
(that enforce the terminal constraint) cannot be employed. However, it
is possible to determine a domain of attraction for the $H_\infty$
controller characterized in \S4.  Consider the following dynamic
programming recursion:
\begin{alignat}{1}
  \label{eq:6.6}
  V_j^0(x) &= \min_{u \in U} \max_{w \in W}\{\ell(x,u,w)+
           V_{j-1}^0(f(x,u,w))\}\\
  \label{eq:6.7}
  \kappa_j(x) &= \arg\min_{u \in U}  \max_{w \in
  W}\{\ell(x,u,w)+V_{j-1}^0(f(x,u,w))\}\\
  \label{eq:6.8}
  X_j^\ast &= \{x \mid f(x,\kappa_j(x),W) \subseteq X_{j-1}^\ast \}
\end{alignat}
with boundary condition
\begin{equation}
  \label{eq:6.9}
  V_0^0(x) = V_f(x), \quad X_0^\ast = X_f
\end{equation}
This is identical to the recursion \eqref{eq:4.4}--\eqref{eq:4.6}
except for the inclusion of the recursion \eqref{eq:6.8} that yields
the sets $X_j^\ast$, $j \ge 0$. Whereas the domain of the value
function $V_j^0\bb$ is $\reals^n$ for all $j \ge 0$, the importance of
the sets $X_j^\ast$ derives from Proposition \ref{p:12} below.
\begin{proposition}
  \label{p:12}
  For every integer $j \ge 0$, every $x \in X_j^\ast$:
  \begin{alignat}{1}
    \label{eq:6.10}
    V_j^0(x) &=V_\infty^0(x)\\
    \label{eq:6.11}
    \kappa_j(x) &=\kappa_\infty(x)
  \end{alignat}
\end{proposition}
\begin{proof}
  Suppose, for some integer $j$, $V_{j-1}^0\bb=V_\infty^0\bb$ on
  $X_{j-1}^\ast$. Then, by \eqref{eq:6.6},
\begin{alignat*}{1}
V_j^0(x) &= \min_{u \in U} \max_{w \in W}\{\ell(x,u,w)+  V_\infty^0(f(x,u,w))\}
=V_\infty^0(x)
\end{alignat*}
for all $x \in X_j^\ast$. Since $V_f(x) =V_\infty^0(x)$ for all $x \in
X_0=X_f$, the desired result follows by induction.
\end{proof}
Hence, the solution to the finite horizon $H_\infty$ problem in \S4
is also the solution to the infinite horizon problem (in the
restricted sets $X_j^\ast$) provided the terminal cost $V_f\bb$
is chosen as described above.  A practical consequence of this result
is that, in computing the value function $V_j^0\bb$ (and
$\kappa_j\bb$), it is only necessary to consider those states lying in
$X_j^\ast \setminus X_{j-1}^\ast$ (since $V_j^0(x)=V_{j-1}^0(x) =
V_\infty(x)$ at all  $x \in X_{j-1}^\ast$).
\begin{definition}
  A set $X$ is robust control invariant for $x^+=f(x,u,w)$ if, for
  every $x \in X$, there exists a $u \in U$ such that $f(x,u,W)
  \subseteq X$.
\end{definition}
It follows from \eqref{eq:6.3} that the set $X_f$ is robust control
invariant.
\begin{theorem}
  \label{th:6}
  (i) The sets $X_j^\ast$ are each robust control invariant and are
  nondecreasing (satisfy $X_j^\ast \subseteq X_{j+1}^\ast$ for all $j
  \ge 0$). (ii) For any $N\ge 0$, the set $X_f$ is finite-time
  attractive with a domain of attraction $X_N^\ast$ for the
  closed-loop system $x^+=f(x,\kappa_\infty(x),w)$. (iii) Suppose
  that, for some finite integer $j \le N$, $X_f$ lies in the interior
  of $X_j^\ast$; then $X_f$ is robustly stable for the system
  $x^+=f(x,\kappa_\infty(x),w)$.
\end{theorem}
\begin{proof}
  (i) Assume $X_{j-1}^\ast$ is robust control invariant. It follows
  from \eqref{eq:6.8} that $X_{j-1}^\ast \subseteq X_j^\ast$ so that
  $X_j^\ast$ is robust control invariant. Since $X_f$ is robust
  control invariant, the desired result follows by induction. (ii) By
  construction, for any integer $j$, and state $x \in X_j^\ast$ is
  robustly steered into $X_{j-1}^\ast$ by the admissible control
  $\kappa_j(x)=\kappa_\infty(x)$. Hence any state $x \in X_N^\ast$ is
  robustly steered into $X_f$ in $N$ steps; the controller $u=K_ux$
  then keeps the state in $X_f$, so that $X_f$ is robustly finite-time
  attractive with a domain of attraction $X_N^\ast$ for the system
  $x^+=f(x,\kappa_\infty(x),w)$. (iii) For any $x \in \reals^n$ let
  $|x|^H \eqbyd d(x,X_f)$ and for any infinite sequence $\{x(i)\}$ in
  $\reals^n$ let $|\{x(i)\}|_\infty^H \eqbyd \sup_{i\ge0} d(x(i),X_f)$.
  From (ii), the controller $\kappa_\infty\bb$ steers any $x \in
  X_j^\ast$ into $X_f$ in no more than $j$ steps and, thereafter,
  keeps the state in $X_f$. Hence, with $x(i) \eqbyd
  \phi(i;x,\pi_\infty,\{w(i)\})$, $\mbf{w} \eqbyd
  \{w(0),w(1),\ldots,w(j-1)\}$, let $\theta:X_j^\ast \times W^j \tends
  \reals$ be defined by
 \[
 \theta(x,\mbf{w}) \eqbyd |x\bb|_\infty^H =
 \max_i\{\phi(i;x,\pi_\infty,\mbf{w}\} \mid i \in
 \{0,1,\ldots,j-1\}\}.
\]
The control law $\kappa_\infty\bb$ is continuous since it is equal to
$\kappa_j\bb$ in $X_j^\ast$; thus $\theta\bb$ is continuous and,
hence, uniformly continuous in $X_j^\ast \times W^j$. Since
$\theta(x,\mbf{w})=0$ for all $x \in X_f \subset X_j^\ast$, all
$\mbf{w} \in W^j$, uniform continuity of $\theta\bb$ implies that, for
all $\varepsilon >0$, there exists a $\delta >0$ such that
$\theta(x,\mbf{w}) < \varepsilon$ ($x(i) \in B_\varepsilon(X_f)$ for
all $i \ge 0$) for  all $(x,\mbf{w}) \in X_j^\ast \times W^j$
satisfying $|x|^H < \delta$ ($x \in B_\varepsilon(X_f)$). This
establishes robust stability of $X_f$.
\end{proof}
The disadvantage of this approach is that the sets $X_j^\ast$ are
obviously subsets of $\reals^n$, the domain of the value functions
$V_j^0\bb$; the sets $X_j^\ast$ are not necessarily convex.

\subsection{$H_\infty$ control: state and control constraints}
Consider the receding horizon controller $u=\kappa_N(x)$. If the
terminal conditions and assumptions stated above in \S6.1 are adopted,
then \\[1ex]
\textbf{C1:} $X_f$ is robust control invariant for $x^+=(A+BK_u)x+Gw$,
$X_f \subseteq X$, $K_u X_f \subseteq U$, $K_w X_f \subseteq W$.\\[0.5ex]
\textbf{C2:} $\min_{u \in U} \max_{w \in W}
\{[\overset{*}{V}_f+\ell](x,u,w) \mid f(x,u,w) \in X_f\} \le 0$ for
all
$x \in X_f$.\\[1ex]
In {C2}, $\overset{*}{V}_f(x,u,w) \eqbyd V_f(f(x,u,w))-V_f(x)$.  If
the recursive dynamic programming equations in \eqref{eq:2.1} --
\eqref{eq:2.3} are employed, we obtain, by a minor modification of the
results in \cite{mayne:2001a},
\S3.3.1, the following results:\\[1ex]
(i:) $X_i$ is robust control invariant for all $i \in \{1,\ldots,
N\}$\\
(ii:) $X_N$ is robust invariant for $x^+=f(x,\kappa_N(x),w)$\\
(iii:) $V_i^0(x) \le V_{i-1}^0(x)\ \forall x \in X_{i-1},\ i \in
\{1,\ldots, N\}$\\
(iv:) $V_N^0(x) \le V_f(x)\ \forall\ x \in X_f$.\\
(v:) The value function satisfies:
\[
[(\sstar{V_N^0}+\ell)\le (V_N^0-V_{N-1}^0)](f(x,\kappa_N(x),w)\le 0
\]
for all $(x,w)\in X_N \times W$. Property (iii) is the
monotonicity property of the value function for the
\emph{constrained}, linear, uncertain system \eqref{eq:1.1} with
cost \eqref{eq:1.5}. Let $V_N(x,\pi,\mbf{w})$ denote the cost if
player $u$ uses the control law $u=\kappa_N(x)$ and the adversary
$w$ uses an arbitrary admissible disturbance sequence $\mbf{w}$.
Then, for any $\ell_2$ disturbance sequence

\begin{multline}
  \label{eq:6.12}
V_N(x,\pi,\mbf{w}) =\sum_{i=0}^{N-1} |y(i)|^2
-(\gamma^2/2) |w(i)|^2  +V_f(x(N))\\ = \sum_{i=0}^{\infty} |y(i)|^2
-(\gamma^2/2) |w(i)|^2 \le V_N^0(x)
\end{multline}
since $\mbf{w}$ is not optimal; here $y(i)=Hz(i)$,
$z(i)=(x(i),u(i))=(x(i),\kappa_N(x(i))$ and $x(i)$ is the solution of
\eqref{eq:1.1} due to initial state $x$, control strategy $\pi$ and
disturbance sequence $\mbf{w}$; we make use of the fact that
$\kappa_N(x) =\kappa_f(x)$ for all $x \in X_f$. It follows that
\begin{equation}
    \label{eq:6.13}
    \sum_{i=0}^\infty |y(i)|^2 \le (\gamma^2/2) \sum_{i=0}^\infty
        |w(i)|^2 + V_N^0(x)
\end{equation}
which is the finite gain property.  Next, if the disturbance is
identically zero,
\[
[(\sstar{V_N^0}+\ell)\le (V_N^0-V_{N-1}^0)](f(x,\kappa_N(x),0))\le 0
\]
so that
\[
V_N^0(f(x,\kappa_N(x),0)) -V_N^0(x) \le -\ell(x,\kappa_N(x),0) \le -(1/2)x'Qx
\]
for all $x \in X_N$ so that the origin is exponentially stable with a
region of attraction $X_N$. Summarizing we have
\begin{theorem}
  \label{th:7}
  The receding horizon controller $u=\kappa_N(x)$ has the following
  properties. The controlled system has the finite $\ell_2$ gain
  property \eqref{eq:6.13} for every initial state in the interior of
  $X_N$ and, if the disturbance is identically zero, the origin is
  exponentially stable with a region of attraction $X_N$.
\end{theorem}
If the disturbance satisfies $w \in W(z)$, where, as before, $z=(x,u)$
and $W$ is such that $w \in W(z)$ implies $|w| \le \delta|z|$, for
some $\delta>0$ (this $W$ models some parametric uncertainties), then
\begin{equation*}
\ell(z,w)=(1/2)|z|_{HH'}^2 -(\gamma^2/2)|w|^2 \ge (c/2)|z|^2 \ge
(c/2)|x|^2
\end{equation*}
for all $z,w$,  some $c>0$, provided that $\delta < (1/\gamma)$.
With this form of bounded disturbance, the origin is robustly,
exponentially stable (the state converges to the origin exponentially
fast despite the disturbance) if, of course,  $\delta < (1/\gamma)$.

We note, in passing, that we can simplify the dynamic programming
recursion, as in \S6.2, by replacing \eqref{eq:2.3} by
\begin{equation*}
   X_j^\ast = X \cap\{x\mid  f(x,\kappa_j(x),W) \subseteq X_{j-1}^\ast\}
\end{equation*}
and the boundary conditions by
\begin{equation*}
  V_0^0(x) = V_f(x), \quad X_0^\ast = X_f
\end{equation*}
However, in this case, \eqref{eq:2.1} remains a constrained
optimization problem because of the state constraint, so the advantage
of using this formulation is not so clear cut. As before,
$X_j^\ast \subseteq X_j$ for each $j$ which introduces
conservatism. However, the sets $X_j^\ast$ are less complex than the
corresponding sets $X_j$.

\section{Conclusion}
We have shown (in \S4) how the solution to the constrained $H_\infty$
problem may be characterized when the system is linear, the cost
quadratic and the constraints polytopic if no state and/or terminad
constraints are present. This characterization required the solution
to a parametric program in which the constraints are polytopic and the
cost piecewise quadratic (rather than quadratic). A novel solution to
this problem is presented in \S3. A characterization of the solution
to the constrained $H_\infty$ problem when state constraints are
present under special (and restrictive) conditions is presented in
\S5; characterization of the solution in the general case appears to
be difficult.  Stability properties of the resultant $H_\infty$
controlled system are briefly discussed in \S6.\\[4ex]

\section*{APPENDICES}

\appendix

\section{Proof of Lemma \ref{p:1}, \S3.1}

We restate Lemma \ref{p:1}:\\[1ex]
\noindent
\textbf{Lemma 1 (Clarke)}
   \emph{
   Suppose $\mc{Z}$ is a polytope in $\reals^n \times \reals^m$ and
   let $\mc{X}$ denote its projection on $\reals^n$ $(\mc{X}= \{x \mid
   \exists u \in \reals^m \textrm{\ such\ that\ } (x,u) \in
   \mc{Z}\})$.  Let $\mc{U}(x) \eqbyd \{u \mid (x,u) \in \mc{Z}\}$.
   Then there exists a $K > 0$ such that, for all $x, x' \in \mc{X}$,
   for all $u \in \mc{U}(x)$,   $d(u,\mc{U}(x')) \le  K|x'-x|$ $($there
   exists a $u' \in \mc{U}(x')$ such that $|u'-u| \le K|x'-x|)$}.

\begin{proof}
  The polytope $\mc{Z}$ is defined by
  \begin{equation*}
   \mc{Z} \eqbyd \{ z= (x,u) \mid Mu \le Nx+p,\  Lx \le l\}.
  \end{equation*}
  where the second set of inequalities is introduced to ensure that no
  row of $M$ is zero ($Lx \le l$ defines the right hand boundary of
  $\mc{Z}$ in Figure \ref{f:1}).  Let $r$ denote the row-dimension of
  $M$, $s$ the row dimension of $L$, let $I \eqbyd \{1,\ldots,r\}$ and
  $J \eqbyd \{1,\ldots,s\}$. For all $(x,u) \in \mc{Z}$, $\mc{U}(x) =
  \{u \mid Mu \le Nx+p\}$ (since $(x,u) \in \mc{Z} \implies Lx \le
  l$). For each $x \in \mc{X}$, $\mc{U}(x)$ is a polytope in
  $\reals^m$. For each $(x,u) \in \mc{Z}$, let $\psi(x,u) \eqbyd \max
  \{M^i u - N^i x- p^i \mid i \in I \}$ and let $\psi^+(x,u) \eqbyd
  \max\{0,\psi(x,u)$ where $M^i$, $N^i$, $p^i$ denote the $i\mathrm{th}$
  row, respectively, of $M$, $N$ and $p$.  Then, for all $x \in
  \mc{X}$, $u \in \mc{U}(x)$ ($(x,u) \in \mc{Z}$) if and only if
  $\psi(z)=0$).  Let $I^0(z) \eqbyd \{i \in I \mid M^i u - N^i x -
  p^i=\psi(z)\}$ index the (most) active constraints.  The associated
  set of (most) active gradients (with respect to $u$) is $\{ g_i \mid
  i \in I^0(z)\}$ where, for each $i$, $g_i= (M^i)'$. Because
  $\mc{U}(x)$ is a polytope, for each $x \in \mc{X}$, the set of
  active gradients is positively linear independent ($ 0 \not \in \co
  \{ g_i(z) \mid i \in I^0(z)\}$ for all $z$ such that $\psi(z)>0$.
  The proximal subgradient of $u \mapsto \psi(x,u)$ is:
  \begin{equation*}
    \delta_p(x,u) \psi(x,u) = \co\{g_i \mid i \in I^0(z)\}.
  \end{equation*}
  and the directional derivative $d\psi(z;h)$ of $\psi\bb$ at $z$ in
  direction $h$ is $\max\{\langle g_i, h\rangle \mid i \in I^0(z)\}$;
  the positive linear independence condition ensures that, at each $z$,
  there exists a direction $h$ along which $\psi(z)$ can be decreased.
  In fact, there exists a $\delta>0$ such that, if $z \in \mc{Z}$ and
  $\psi(z) > 0$ ($z \not \in \mc{Z}$), then $\zeta \in \delta_p(z)$
  implies $|\zeta| \ge \delta$.  So, by Theorem 3.1 in
  \cite{clarke:ledyaev:stern:wolenski:1998}, for all $(x,u) \in
  \mc{Z}$, all $x' \in \mc{X}$, $d(u,\mc{U}(x')) \le \psi(x',u)/\delta
  \le (c/\delta) |x-x'|$ ($c \eqbyd \max \{|N^j| \mid j \in J\}$)
  since $\psi(x',u) \le \psi(x,u) + \max \{N^j(x-x')\mid j \in J\} \le c
  |x-x'|\}$.  This proves the lemma with $K \eqbyd c/\delta$.
\end{proof}

\section{Continuous differentiability of the value function of
  a parametric quadratic program}

See \cite{mayne:2001a} and \cite{borrelli:2002} for related results.
We consider the standard parametric quadratic program:
\begin{equation}
  \label{eq:A.1}
  \mbb{P}(x): \quad V^0(x)=\min_{u} \{V(x,{u}) \mid (x,u) \in
  \mc{Z}\}.
\end{equation}
where $x \in \reals^n$, ${u} \in \reals^l$ and $\mc{Z}$ is a
polytope with a non-empty interior.
We assume\\[1ex]
\textbf{A1}: The cost function $V\bb$ is strictly convex and
continuously differentiable.
The constraint $(x,u) \in \mc{Z}$ imposes an implicit constraint
on ${u} \in \mc{U}(x)$  where the set-valued function $\mc{U}\bb$
is defined by

\begin{equation}
  \label{eq:A.2}
  \mc{U}(x) = \{{u} \mid (x,u) \in \mc{Z}\} = \{{u} \mid
  M{u} \le Nx+p\}
\end{equation}
so that $\mbb{P}(x)$ may be written in the form
\begin{equation}
  \label{eq:A.3}
   \mbb{P}(x): \quad V^0(x)=\min_{u} \{V(x,{u}) \mid {u}  \in
   \mc{U}(x)\}
\end{equation}
The (unique) solution of $\mbb{P}(x)$, for each $x \in \mc{X}$, is
\begin{equation}
  \label{eq:A.4}
  {u}^0(x) = \arg\min_{u} \{V(x,{u}) \mid {u}  \in
   \mc{U}(x)\}
\end{equation}
The domain of $V^0\bb$, ${u}\bb$ and $\mc{U}\bb$,  is the polytope
\begin{equation}
   \label{eq:A.5}
  \mc{X} =\{x \mid \mc{U}(x) \not = \emptyset\} = \{x \mid \exists\ {u}
  \in \mc{U}(x) \} =\proj_X (\mc{Z})
\end{equation}
Let $p\ge l$ denote the number of rows of $M,N$ and $c$.  It is known
that $V^0\bb$ is continuous piecewise quadratic and continuous and
${u}^0\bb$ is piecewise affine and continuous, being quadratic and
affine, respectively, in the polytopes $R_I$, $I \subseteq \mc{I}_p$ that
constitute a polytopic partition of $\mc{X}$. Each region is
characterized by a set of active constraints $I$, i.e. for all $x \in
R_I$:
\begin{alignat}{1}
    \label{eq:A.6}
    M_I{u}^0(x)&=N_Ix+ c_I  \\
    \label{eq:A.7}
    M_i{u}&\le N_ix+ p_i \mathrm{\ for\ all\ } i \in I^c\\
    \label{eq:A.8}
    -\nabla_u V(x,{u}^0(x)) &\in  PC_I(x)
\end{alignat}
where $M_I$, $N_I$ and $c_I$ denote the matrices with, rows $M_i$,
$N_i$ and $c_i$, respectively, $i \in I$, and $PC_I(x) \eqbyd
\{M_I'\lambda \mid \lambda \ge 0\}$ is the polar cone to the cone $F(x)
\eqbyd \{h \mid M_I h \le 0\}$ of feasible directions at $x$; for each
$I \subseteq \mc{I}_p$, $I^c$ denotes the complement of $I$ in $\mc{I}_p$.
Thus $V^0\bb$ is continuously differentiable (in fact analytic) in the
interior of each region $R_I$, $I \subseteq \mc{I}_p$. We may assume,
without loss of generality, that $M_I$ has maximal rank.
Our final assumption is:\\[1ex]
\textbf{A2:} For any two adjacent regions $R_{I_1}$ and $R_{I_2}$
($R_{I_1} \cap R_{I_2} \not = \emptyset$) either ${I_1} \subset
{I_2}$ or ${I_1} \supset {I_2}$.\\[1ex]
Assumption \textbf{A2} will often be satisfied, but there do exist
counterexamples.
\begin{theorem}
  \label{th:A1}
  Suppose $V\bb$ is continuously differentiable and that assumptions
  \textbf{A1} -- \textbf{A2} are satisfied.  Then $V^0\bb$ is
  continuously differentiable in $\mc{X}$.
\end{theorem}
\begin{proof}
  It is known that $V^0\bb$ is continuous piecewise quadratic and
  continuous and ${u}^0\bb$ is piecewise affine and continuous, being
  quadratic and affine, respectively, in the polytopes $R_I$, $I
  \subseteq \mc{I}_p$ that constitute a polytopic partition of
  $\mc{X}$.  Each region is characterized by a set of active
  constraints $I$, i.e.  $R_I$ is defined by the inequalities
  \eqref{eq:A.6}-\eqref{eq:A.8}.  Thus $V^0\bb$ is continuously
  differentiable (in fact analytic) in the interior of each region
  $R_I$, $I \subseteq \mc{I}_p$.  Consider the continuous
  differentiability of $V^0\bb$ on the boundary between two regions
  $R_{I_1}$ and $R_{I_2}$ say where $I_1 \subseteq I_2$. For any $I
  \subseteq \mc{I}_p$ such that $R_I \neq \emptyset$, any $x \in
  R_{I}$,
\begin{equation}
    \label{eq:A.9}
  {u}^0(x)= \tmbf{u}_{I}^0(x)+\bmbf{u}_{I}^0(x)
\end{equation}
where, for each index set $I$, $\tmbf{u}_{I}^0(x) \in
\mathrm{range}(M_{I}')$$=\mc{N}(M_I)^\perp$ (the row space of $M_I$)
is that $u$ of minimum norm satisfying $M_{I}{u}=N_{I}x+ c_{I}$
and $\bmbf{u}_{I}^0(x) \in \mathrm{range}(M_{I}')^\perp=\mc{N}(M_I)$
($\mc{N}(M_I)$ is the null space of $M_I$).  Roughly speaking,
$\tmbf{u}_{I}^0(x)$ satisfies the constraints, and $\bmbf{u}_{I}^0(x)$
optimizes. It is easily shown that both $\tmbf{u}_{I}\bb$ and
$\bmbf{u}_{I}^0\bb£$ are affine in $x$, satisfying, respectively
\begin{equation}
  \label{eq:a.10}
 \tmbf{u}_{I}^0(x) = \tilde{K}_{I}x +\tilde{k}_{I},\quad
 \bmbf{u}_{I}^0(x) = \bar{K}_{I}x +\bar{k}_{I}
\end{equation}
where $\tilde{k}_{I}$ and the columns of $\tilde{K}_{I}$ lie in
$\mathrm{range}(M_{I}')$ and $\bar{k}_{I} $ and the columns of
$\bar{K}_{I}$ lie in $\mathrm{range}(M_{I}')^\perp$.
In fact, $\tilde{K}_{I}$ and $\tilde{k}_{I}$ are given by
\begin{equation}
  \label{eq:A.11}
  \tilde{K}_{I} = M_{I}^\dag N_{I},\quad \tilde{k}_{I}
    =M_{I}^\dag c_{I}
\end{equation}
where $M_I^\dag$, the Moore-Penrose pseudo inverse of $M_I$, is given
by
\begin{equation}
  \label{eq:A.12}
  M_I^\dag = (M_I'M_I)^{-1}M_I'
\end{equation}
Since $V^0\bb$ is continuously differentiable in each region $R_I$,
consider the continuous differentiability of $V^0\bb$ on the boundary
between two regions, $R_{I_1}$ and $R_{I_2}$ say, where $I_1 \subseteq
I_2$ Because, for $x \in R_{I}$, ${u}^0(x)$ minimizes (with
respect to ${u}$) the continuously differentiable function
$V(x,{u})$ in the hyperplane $\{{u} \mid M_{I_1}{u} =
N_{I_1}x+c_{I_1}\}$, we have
\begin{equation}
  \label{eq:A.13}
  -\nabla_{u}V(x,{u}^0(x)) \in \{M_{I_1}'\lambda \mid \lambda \ge 0\}
  \subseteq \mathrm{range}(M_{I_1}')
\end{equation}
Hence
\begin{equation}
  \label{eq:A.14}
  (\partial/\partial x)V_x^0(x)=  (\partial/\partial x)V(x,{u}^0(x))+
  (\partial/\partial u)V(x,{u}^0(x))(\partial/\partial
  x)\tmbf{u}_{I_1}^0(x)
\end{equation}
since $(\partial/\partial u)V((x,{u}^0(x)))(\partial/\partial
x)\bmbf{u}_{I_1}^0(x)=0$ (because of \eqref{eq:A.13} and the fact that
$(\partial/\partial x)\bmbf{u}_{I_1}^0(x)$ $=\bar{K}_{I_1}$) and the
columns of $\bar{K}_{I_1}$ lie in $\mathrm{range}(M_{I_1}')^\perp$.
Suppose now $x \tends x^\ast \in R_{I_1} \cap R_{I_2}$, $x \in
R_{I_1}$. Then
 \begin{equation}
   \label{eq:A.15}
   (\partial/\partial x)V^0(x) \tends (\partial/\partial x)
   V(x^\ast,{u}^0(x^\ast))+ (\partial/\partial u)
   V(x,{u}^0(x^\ast))\tilde{K}_{I_1}
 \end{equation}
 where both $\nabla_u V(x,{u}^0(x^\ast))'$ and the columns of
 $\tilde{K}_{I_1}=(\partial/\partial x)\tmbf{u}_{I_1}^0(x^\ast)$ lie in
 $\mathrm{range}(M_{I_1}')$.  Next consider a $x \in R_{I_2}$ such
 that $x \tends x^\ast \in R_{I_1} \cap R_{I_2}$. Arguing as above we
 deduce
 \begin{equation}
   \label{eq:A.16}
   (\partial/\partial x)V_x^0(x) \tends (\partial/\partial x)
   V(x^\ast,{u}^0(x^\ast))+
  (\partial/\partial u)V_u((x,{u}^0(x^\ast))\tilde{K}_{I_2}
 \end{equation}
 where $\nabla_uV(x,{u}^0(x^\ast))$ lies in
 $\mathrm{range}(M_{I_1}')$ (as above) but the columns of
 $\bar{K}_{I_2}=$ $(\partial/\partial x)\tmbf{u}_{I_2}^0(x^\ast)$ lie
 in $\mathrm{range}(M_{I_2}')$.  We show below that $I_1 \subseteq
 I_2$ implies that $M_{I_1}\tilde{K}_{I_1} = M_{I_1}\tilde{K}_{I_2}$.
 Since $\nabla_uV((x,{u}^0(x^\ast))$ lies in
 $\mathrm{range}(M_{I_1}')$, it follows that
 \begin{equation}
  \label{eq:A.17}
  (\partial/\partial u) V(x,{u}^0(x^\ast))\tilde{K}_{I_1}
   = (\partial/\partial u)V_u((x,{u}^0(x^\ast))\tilde{K}_{I_2}
 \end{equation}
Equations \eqref{eq:A.15} - \eqref{eq:A.17} establish the continuous
differentiability of $V^0\bb$ at $x^\ast \in R_{I_1} \cap R_{I_2}$.

We have now to show that $I_1 \subseteq I_2$ implies that
$M_{I_1}\tilde{K}_{I_1} = M_{I_1}\tilde{K}_{I_2}$.
Suppose
\begin{equation*}
 M_{I_2}= \left[ \begin{array}{l} M_{I_1} \\m \end{array} \right],
 \quad N_{I_2}= \left[ \begin{array}{l} N_{I_1} \\n \end{array} \right],
\end{equation*}
Then, from \eqref{eq:A.11}
\begin{alignat*}{1}
  M_{I_1} \tilde{K}_{I_2} & = M_{I_1} M_{I_2}^\dag N_{I_2}
\end{alignat*}
But
\begin{equation*}
     M_{I_2} \tilde{K}_{I_2}=N_{I_2}
\end{equation*}
so that
\begin{equation*}
  \left[ \begin{array}{c} M_{I_1}\\m \end{array} \right]
  \tilde{K}_{I_2} = \left[ \begin{array}{c} N_{I_1}\\n \end{array} \right]
\end{equation*}
from which it follows that
\begin{equation*}
  M_{I_1}\tilde{K}_{I_2} =N_{I_1}= M_{I_1}\tilde{K}_{I_1}.
\end{equation*}
It follows from \textbf{A2} that $V^0\bb$ is continuously
differentiable in $\mc{X}$.

\end{proof}

\end{document}